\documentclass{amsart}%
\usepackage{palatino, mathpazo}
\usepackage{amsfonts}
\usepackage{amsmath}
\usepackage{amssymb,latexsym,xcolor}
\usepackage{graphicx}
\usepackage[mathscr]{eucal}
\usepackage{harpoon}
\usepackage{amssymb}
\usepackage[
linktocpage=true,colorlinks,citecolor=magenta,linkcolor=blue,urlcolor=magenta]{hyperref}
\allowdisplaybreaks[4]
\providecommand{\U}[1]{\protect \rule{.1in}{.1in}}

\makeatletter

\def\Rmnum#1{\expandafter\@slowromancap\romannumeral #1@}
\makeatother

\newtheorem{theorem}{Theorem}[section]

\newtheorem{corollary}[theorem]{Corollary}
\newtheorem{definition}[theorem]{Definition}
\newtheorem{example}[theorem]{Example}
\newtheorem{Lemma}[theorem]{Lemma}
\newtheorem{Proposition}[theorem]{Proposition}

\theoremstyle{remark}
\newtheorem{remark}[theorem]{Remark}

\numberwithin{equation}{section}

\def\C{\mathbb C}

\def\N{\mathbb N}

\def\R{\mathbb R}

\def\al{\alpha}
\def\ga{\gamma}

\def\e{\varepsilon}

\def\iy{\infty}

\def\la{\lambda}

\def\vp{\varphi}
\def\pa{\partial}

\def\nm#1{\|#1\|}

\def\lab{\label}
\def\f{\frac}
\def\bt{\begin{theorem}}
\def\et{\end{theorem}}
\def\bl{\begin{Lemma}}
\def\el{\end{Lemma}}
\def\bd{\begin{definition}}
\def\ed{\end{definition}}
\def\bc{\begin{corollary}}
\def\ec{\end{corollary}}
\def\bpr{\begin{Proposition}}
\def\epr{\end{Proposition}}
\def\bp{\begin{proof}}
\def\ep{\end{proof}}
\def\bx{\begin{example}}
\def\ex{\end{example}}
\def\br{\begin{remark}}
\def\er{\end{remark}}
\def\be{\begin{equation}}
\def\ee{\end{equation}}
\def\bal{\begin{align}}
\def\bn{\begin{enumerate}}
\def\en{\end{enumerate}}
\def\eal{\end{align}}
\def\bg{\begin{align*}}
\def\eg{\end{align*}}
\def\bcs{\begin{cases}}
\def\ecs{\end{cases}}
\def\abs#1{\lvert#1\rvert}
\def\RNum#1{\uppercase\expandafter{\romannumeral #1\relax}}

\def\CR{\mathcal C}

\def\MR{{\mathcal M}}

\def\PR{{\mathcal P}}
\def\SR{{\mathcal S}}
\def\TR{{\mathcal T}}

\def\ZR{{\mathcal Z}}

\def\bean{\begin{eqnarray*}}
\def\eean{\end{eqnarray*}}
\def\loc{\operatorname{\rm loc}}

\def\N{{\mathbb N}}

\def\C{{\mathbb C}}
\def\bean{\begin{eqnarray*}}
\def\eean{\end{eqnarray*}}
\def\loc{\operatorname{\rm loc}}

\def\sbr#1{\left(#1\right)}
\def\mbr#1{\left[#1\right]}
\def\lbr#1{\left\{#1\right\}}
\def\abr#1{\left\langle#1\right\rangle}
\def\rd{{\mathrm d}}
\def\abs#1{\left\lvert#1\right\rvert}

\def\wt{\widetilde}
\def\lab{\label}

\begin{document}
\title[Lane-Emden equation with vanishing potentials]{Energy quantization of the two dimensional Lane-Emden equation with vanishing potentials}
\author{Zhijie Chen}
\address{Department of Mathematical Sciences, Yau Mathematical Sciences Center,
Tsinghua University, Beijing, 100084, China}
\email{zjchen2016@tsinghua.edu.cn}
\author{Houwang Li}
\address{Yanqi Lake Beijing Institute of Mathematical Sciences and Applications, Beijing, 101408, China}
\email{lhwmath@bimsa.cn}



\begin{abstract}
We study the concentration phenomenon of the Lane-Emden equation with vanishing potentials
	$$\bcs
		-\Delta u_n=W_n(x)u_n^{p_n},\quad u_n>0,\quad\text{in}~\Omega,\\
		u_n=0,\quad\text{on}~\pa\Omega,\\
		\int_\Omega p_n W_n(x)u_n^{p_n}\rd x\le C,
	\ecs$$
where $\Omega$ is a smooth bounded domain in $\mathbb{R}^2$, $W_n(x)\geq 0$ are bounded functions with zeros in $\Omega$, and $p_n\to\iy$ as $n\to\iy$. A typical example is $W_n(x)=|x|^{2\al}$ with $0\in\Omega$, i.e. the equation turns to be the well-known H\'enon equation. The asymptotic behavior for $\al=0$ has been well studied in the literature. While for $\al>0$, the problem becomes much more complicated since a singular Liouville equation appears as a limit problem. In this paper, we study the case $\al>0$ and prove a quantization property (suppose $0$ is a concentration point)
	\[p_n|x|^{2\alpha}u_n(x)^{p_n-1+t}\to 8\pi e^{\frac{t}{2}}\sum_{i=1}^k\delta_{a_i}+8\pi(1+\al)e^{\f{t}{2}}c^t\delta_0, \quad t=0,1,2,\]
for some $k\ge0$, $a_i\in\Omega\setminus\{0\}$
and some $c\ge1$. Moreover, for $\al\not\in\N$, we show that the blow up must be simple, i.e. $c=1$. As applications, we also obtain the complete asymptotic behavior of ground state solutions for the H\'enon equation.
\end{abstract}
\maketitle

\section{Introduction}
\label{section-1}

In the past years much attention has been paid to the blow-up analysis for solution sequences $u_n(x)$ of the Lane-Emden type equation
\be\lab{equ-1} \begin{cases}
	-\Delta u_n=W_n(x)|u_n|^{p_n-1}u_n,\quad\text{in}~\Omega,\\
	u_n=0,\quad\text{on}~\pa\Omega,\\
	\int_\Omega p_n W_n(x)|u_n|^{p_n}\rd x\le C,
\end{cases}\ee
where $\Omega\subset\R^2$ is a smooth bounded domain, $W_n(x)\geq 0$ are bounded functions with zeros in $\Omega$, $p_n\to\iy$ as $n\to\iy$, and $C>0$ is a constant independent of $n$. As in the literature,  $a\in\overline\Omega$ is called a blow-up point of $p_nu_n$ if there exists $\lbr{x_n}\subset\Omega$ such that $x_n\to a$ and $p_nu_n(x_n)\to\iy$. In this case, we also call this $a$ a blow-up point of $u_n$ for convenience.

When $W_n(x)\equiv1$, \eqref{equ-1} turns to be the well known two demensional Lane-Emden equation
\be\lab{equ-2} \begin{cases}
	-\Delta u_n=|u_n|^{p_n-1}u_n,\quad\text{in}~\Omega,\\
	u_n=0,\quad\text{on}~\pa\Omega,\\
	\int_\Omega p_n|u_n|^{p_n}\rd x\le C.
\end{cases}\ee
The asymptotic behaviors of \emph{positive solutions} of \eqref{equ-2} have been well studied by various mathematicians in a series of papers \cite{LE-8,LE-1,asy2-4,LE-5,LE-6}, and the main results can be summarized as follows: Let $u_n$ be a sequence of positive solutions of \eqref{equ-2}. Then there exists a finite set $\SR=\lbr{a_1,\cdots,a_k}\subset\Omega$ consisting of blow-up points of $p_nu_n$ such that up to a subsequence, for a suitable $r_0>0$,
\be\label{1-3}
	\sup_{B_{r_0}(a_i)} u_n(x)\to\sqrt{e}\quad\text{and}\quad \sup_{\overline{\Omega}\setminus \bigcup_{i=1}^k B_r(a_i)}p_nu_n(x)\le C_r, \quad\text{\it for any}~r>0,
\ee
\be\label{1-4}
	p_nu_n(x)^{p_n-1}\to 8\pi\sum_{i=1}^k \delta_{a_i}\quad\text{\it weakly in the sense of measures},
\ee
where $\delta_{a_i}$ is the Dirac measure at $a_i$, and \[B_r(a):=\{x\in\mathbb{R}^2\,:\, |x-a|<r\},\quad B_r:=B_r(0)\] denote open balls.
Furthermore, for each $1\leq i\leq k$, a suitable scaling of $u_n$ near $a_i$ converges in $\CR_{loc}^2(\R^2)$ to an entire solution $U$ of the Liouville equation
\be\lab{Liouville-1} \begin{cases}
	-\Delta U=e^U\quad\text{in}~\R^2,\\
	\int_{\R^2}e^U\rd x<\iy.
\end{cases}\ee

On the other hand, the asymptotic behaviors of \emph{nodal solutions} of \eqref{equ-2} are much more difficult to study and there are only some partial results; see \cite{LE-7,LE-4,LE-11}. In particular, comparing to positive solutions, new phenomena appear for nodal solutions. For example, Grossi, Grumiau and Pacella \cite{LE-11} studied \emph{the least energy radial nodal solutions} in a ball, and proved that the limit profile of these nodal solutions looks like a superposition of two bubbles, one related to a regular limit problem \eqref{Liouville-1} and another one related to a singular limit problem \be\lab{Liouville-2} \begin{cases}
	-\Delta U=e^U+H\delta_0,\quad\text{in}~\R^2,\\
	\int_{\R^2}e^U\rd x<\iy,
\end{cases}\ee
where $H$ is a suitable constant. More precisely, a suitable scaling of the positive parts $u_n^+=\max\{u_n, 0\}$ converges to a solution of the Liouville equation \eqref{Liouville-1}, while a suitable scaling of the negative parts $u_n^-=\min\{u_n,0\}$ converges to a singular solution of the singular Liouville equation \eqref{Liouville-2}.

One purpose of this paper is to show that that for \emph{positive solutions} of \eqref{equ-1}, if $W_n(x)$ vanishes (with finite order) at some points, then the singular Liouville equation \eqref{Liouville-2} appears again as a limit problem. This is a different feature comparing to positive solutions of the Lane-Emden equation \eqref{equ-2}.

Our another interest of studying \eqref{equ-1} is originated from the  H\'enon equation
\be \begin{cases}
	-\Delta u_n=|x|^{2\al}|u_n|^{p_n-1}u_n,\quad\text{in}~B_1,\\
	u_n=0,\quad\text{on}~\pa B_1,
\end{cases} \ee
which was introduced by H\'enon \cite{Henon-0} in the study of stellar clusters in radially symmetric settings in 1973.
Here we consider more general potentials $W_n(x)$. Suppose $W_n(x)$ has the form
\be\lab{con-1}      W_n(x)=\overline W_{n}(x)\prod_{i=1}^{m}|x-q_i|^{2\al_i}, \ee
where $m\ge1$, $\al_i>0$ and $\overline W_n$ satisfies 
\be\lab{con-2}  	0<\f{1}{C}\le \overline W_{n}(x)\le C<\iy,\quad |\nabla \overline W_{n}(x)|\le C,\quad\text{for}~x\in\Omega,  \ee
for some positive  constant $C$ independent of $n$. Denote the zero set of $W_n(x)$ by
\be\lab{Z} \ZR:=\lbr{x\in\Omega:~W_n(x)=0}=\lbr{q_1,\cdots,q_m}. \ee
We will see that the problem will become very subtle if $p_nu_n$ blows up at some points in $\ZR$.

\subsection{Local problems}

We start from a Br\'ezis-Merle type result. In \cite{MF-1} Br\'ezis-Merle gave their famous alternative results for the Liouville problem $-\Delta u_n=V_n(x)e^{u_n}$ in $\Omega$. Later, Ren-Wei \cite{LE-5,LE-6} developed their method to handle the least energy solutions of the Lane-Emden equation \eqref{equ-2}. Here we follow Ren-Wei's idea to prove the following Br\'ezis-Merle type result.
\bt\lab{thm-BM}
Suppose $p_n\to\iy$ and $u_n$ is a solution sequence of 
\be\lab{equ-BM} \begin{cases}
	-\Delta u_n=V_n(x)u_n^{p_n},\quad u_n>0,\quad\text{in}~\Omega,\\
	\int_\Omega p_n V_n(x)u_n^{p_n}\rd x\le C.
\end{cases}\ee
Then under the condition that
\be\lab{con-BM} 0\le V_n(x)\le C,\quad |\nabla V_n(x)|\le C,\quad\text{for all}~x\in\Omega, \ee
after passing to a subsequence (still called $u_n$), one of the following alternatives holds:
\begin{itemize}
\item[(i)] $u_n \to 0$ uniformly in $L_{loc}^\iy(\Omega)$ with $\nm{p_nu_n}_{L^\iy(K)}\leq C_K$ for any compact subset $K\Subset\Omega$.
\item[(ii)] There exist a non-empty finite set $\Sigma=\lbr{a_1,\cdots,a_k}\subset\Omega$ and corresponding sequences $\lbr{x_{n,i}}_{n\in\N}$ in $\Omega$ for $i=1,\cdots,k$, such that $x_{n,i}\to a_i$ and $u_n(x_{n,i})\to\ga_i\ge 1$ as $n\to\iy$. Moreover, $\nm{p_nu_n}_{L^\iy(K)}\le C_K$ for any compact subset $K\Subset\Omega\setminus\Sigma$, and 
	\be\label{lambda-i} p_nV_n(x)u_n(x)^{p_n-1}\to\sum_{i=1}^k\beta_{a_i}\delta_{a_i},\quad p_nV_n(x)u_n(x)^{p_n}\to\sum_{i=1}^k\la_{a_i}\delta_{a_i}\ee
weakly in the sense of measures in $\Omega$ with $\beta_{a_i}\ge\f{4\pi e}{\ga_i}$ and $\la_{a_i}\ge4\pi e$.
\end{itemize}
\et

Note that we need no boundary conditions on $u_n$ in Theorem \ref{thm-BM}. When the alternative (ii) holds, the set $\Sigma$ only consists of those blow-up points of $p_nu_n$ contained in $\Omega$, i.e. whether $p_nu_n$ blows up at some points of $\partial\Omega$ or not is unknown. 

After Theorem \ref{thm-BM}, a natural question arises:

\vskip0.1in
\noindent{\bf Question.} \emph{When the alternative (ii) holds, can one compute the exact values of $\beta_{a_i}$, $\la_{a_i}$ for every $i$}? 
\vskip0.1in

An easy situation is that $V_n(x)$ is bounded below away from zero near the blow-up point $a_i$, and we call this \emph{a regular case}.

\bt\lab{thm-regular}
Suppose $p_n\to\iy$, $r>0$ and $u_n$ is a solution sequence of 
\be\lab{equ-regular} \begin{cases}
	-\Delta u_n=V_n(x)u_n^{p_n},\quad u_n>0,\quad\text{in}~ B_r,\\
	\int_{B_r} p_n V_n(x)u_n^{p_n}\rd x\le C,
\end{cases}\ee
with $0$ being the only blow-up point of $p_nu_n$ in $B_r$, i.e.,
\be\lab{con-regular1} 	\max_{B_r}p_nu_n\to\iy  \quad \text{and}\quad \max_{\overline B_r\setminus B_\delta} p_nu_n\le C_\delta,\quad\text{for any}~0<\delta<r. \ee
Then under the condition that
\be\lab{con-regular2} 0<\f{1}{C}\le V_n(x)\le C,\quad |\nabla V_n(x)|\le C,\quad\text{for}~x\in B_r, \ee
after passing to a subsequence (still called $u_n$), it hold $\max_{B_r}u_n\to\sqrt e$ and
\be\lab{result-regular}  p_n|x|^{2\alpha}V_n(x)u_n(x)^{p_n-1+t}\to 8\pi e^{\f{t}{2}}\delta_0, \quad t=0,1,2  \ee
weakly in the sense of measures.
\et

Theorem \ref{thm-regular} improves \cite[Theorem 1.1]{LE-1} in the sense that $V_n\not\equiv1$ is allowed and no boundary condition $u_n=0$ is needed. The idea of proving Theorem \ref{thm-regular} is similar to that of \cite[Theorem 1.1]{LE-1}, and for the reader's convenience we will give the proof in Section 3. First, by the blow-up analysis around a local maximum of $u_n$, we are led to a solution $U_0$ of the Liouville equation \eqref{Liouville-1}. The classical result of Chen-Li \cite{classification-1} characterizes all these solutions, which implies 
$\int_{\R^2}e^{U_0}=8\pi$. Then, by the local Pohozaev identity and the Green's representation formula, we get a decay estimate of $u_n$, which is used to apply the Dominate Convergence Theorem to get the convergence of energies, and hence get the desired results. 

\vskip0.1in
Now a delicate situation is that $V_n(x)$ vanishes (with finite order) at a blow-up point $a_i$, and we call this \emph{a singular case} since the blow-up around a local maximum of $u_n$ near $a_i$ will lead to the singular Liouville problem. Due to $\al\in\N$ or not, we have different results.

\bt\lab{thm-singular}
Suppose $p_n\to\iy$, $r>0$, $\al>0$ and $u_n$ is a solution sequence of 
\be\lab{equ-singular} \begin{cases}
	-\Delta u_n=|x|^{2\al}V_n(x)u_n^{p_n},\quad u_n>0,\quad\text{in}~ B_r,\\
	\int_{B_r} p_n |x|^{2\al}V_n(x)u_n^{p_n}\rd x\le C,
\end{cases}\ee
with $0$ being the only blow-up point of $p_nu_n$ in $B_r$, i.e.,
\be\lab{con-singular1} 	\max_{B_r}p_nu_n\to\iy  \quad \text{and}\quad \max_{\overline B_r\setminus B_\delta} p_nu_n\le C_\delta,\quad\text{for any}~0<\delta<r. \ee
Then under the condition that
\be\lab{con-singular2} 0<\f{1}{C}\le V_n(x)\le C,\quad |\nabla V_n(x)|\le C,\quad\text{for}~x\in B_r, \ee
after passing to a subsequence (still called $u_n$), it holds that $\max_{B_r}u_n\to\ga\ge\sqrt e$ and
\be\lab{result-singular1}  p_n|x|^{2\alpha}V_n(x)u_n(x)^{p_n-1+t}\to 8\pi(1+\al)e^{\f{t}{2}}c^t\delta_0, \quad t=0,1,2,\ee
weakly in the sense of measures for some $c\in[1,\ga]$. Moreover, there holds $c=1$ and $\ga=\sqrt e$ if $\al\not\in\N$.
\et

\begin{remark} 
It is interesting to compare Theorem \ref{thm-singular} with the results for singular mean field problems. Suppose $u_n$ solves
	\[-\Delta u_n=|x|^{2\al}V_ne^{u_n}\quad\text{\it in }\;B_r,\] 
and assume that $0$ is the only blow-up point. Then under the condition \eqref{con-singular2}, Tarantello \cite{SMF-1} proved that $|x|^{2\al}V_ne^{u_n}\to\beta_0\delta_0$ with $\beta_0\in 8\pi\N_{\geq 1}\cup\{8\pi(1+\al)+8\pi\N\}$. There are also explicit examples in \cite{SMF-1} to show that $\beta_0$ can take any value contained in $8\pi\N_{\geq 1}\cup\{8\pi(1+\al)+8\pi\N\}$. Here we get a quite surprising result, that is the energy must be $8\pi(1+\al)$, i.e.,
	$$p_n|x|^{2\alpha}V_n(x)u_n(x)^{p_n-1}\to 8\pi(1+\al)\delta_0,$$
for any $\al>0$. 
\end{remark}

The proof of Theorem \ref{thm-singular} is much more complicated than that of Theorem \ref{thm-regular}, and the main difficulty is the lack of the following condition 
\be\lab{result-singular2} 
	\sup_{B_{\e_0}}p_n|x|^{2+2\al}V_n(x)u_n(x)^{p_n-1}\le C, \quad\text{for some}~\e_0\in(0,r).
\ee

The proof consists of two main ingredients. First, we assume \eqref{result-singular2} holds, and then by a blow-up around a local maximum of $u_n$, we are led to a solution $U_\al$ of the singular Liouville problem
\be\lab{Liouville-3}\begin{cases}
	-\Delta U_\al=|x|^{2\al}e^{U_\al},\quad\text{in}~\R^2,\\
	\int_{\R^2}|x|^{2\al}e^{U_\al}\rd x<\iy.
\end{cases}\ee
Since $\Delta(\f{1}{2\pi}\ln|x|)=\delta_0$, we see that \eqref{Liouville-2} and \eqref{Liouville-3} are equivalent in the sense that $U_\al$ is a solution of \eqref{Liouville-3} if and only if $U_\al+2\al\ln|x|$ is a solution of \eqref{Liouville-2} with $H=-4\pi\al$. A result of Prajapat-Tarantello \cite{classification-2} (see \cite{classification-1} for $\al=0$) characterizes all solutions of  \eqref{Liouville-3}, from which we know that
\be  \int_{\R^2}|x|^{2\al}e^{U_\al}\rd x=8\pi(1+\al). \ee
Then by the local Pohozaev identity and the Green's representation formula, we get a decay estimate of $u_n$, and hence we obtain $\beta_t=8\pi(1+\al)e^{\frac{t}2}$ for $t=0,1,2$. 

Second, we assume \eqref{result-singular2} does not hold. Note that equation \eqref{equ-singular} is formally invariant under the transformation
\be\lab{tem-101} v_n(x)=r^{\al_n}u_n(rx),\quad\text{with}~\al_n=\f{2+2\al}{p_n-1}. \ee
Thanks to this transformation and inspired by \cite{SMF-1}, we can construct a decomposition of $u_n$; see Proposition \ref{decomposition}. In this direction, we reduce the singular case to some regular cases. By accurate analysis, we show that there is no energy loss in neck domains. Then using Theorems \ref{thm-BM} and \ref{thm-regular}, we compute the exact values of the correponding energies, which gives 
	$$p_n|x|^{2\alpha}V_n(x)u_n(x)^{p_n-1+t}\to 8\pi e^{\f{t}{2}}(\sum_{i=1}^l N_ic_i^t)\delta_0,\quad t=0,1,2,$$
for some $l\ge1$, $c_i\ge1$, $N_i\in\N$ for $i=1,\cdots,l$. Then comparing these energies by Pohozaev identity, we get $l=1$ and $N_1=1+\al$. So that if $\al\not\in\N$, we get a contradiction, and then condition \eqref{result-singular2} holds, which gives Theorem \ref{thm-singular} for $\al\not\in\N$. While for $\al\in\N$, the result is more complicated. For the mean field equation with integer singular sources $\al\in\N$, Kuo-Lin \cite{SMF-4} and Bartolucci-Tarantello \cite{SMF-5} showed the non-simple blow-up  phenomena happens, i.e., condition \eqref{result-singular2} does not hold. We refer to \cite{Nonsimple-1,Nonsimple-2} for more information of the non-simple blow-up. Hence it is an interesting problem to consider the blow-up phenomena of \eqref{equ-singular} with $\al\in\N$, and in a following paper, we would like to study this case.

\subsection{Boundary value problems}
Thanks to the above local properties, we are in position to study positive solutions of our initial problem \eqref{equ-1}, i.e.
\be\lab{equ-1-0} \begin{cases}
	-\Delta u_n=W_n(x)u_n^{p_n},\quad u_n>0,\quad\text{in}~\Omega,\\
	u_n=0,\quad\text{on}~\pa\Omega,\\
	\int_\Omega p_n W_n(x)u_n^{p_n}\rd x\le C,
\end{cases}\ee
with $W_n(x)$ satisfying \eqref{con-1}-\eqref{con-2}.
For a solution sequence $u_n$ of \eqref{equ-1-0}, we define the set $\SR$ of blow-up points of $p_nu_n$ as
\be\lab{S}  \SR:=\lbr{a\in\overline\Omega:~\exists \{x_n\}\subset\Omega,~x_n\to a,~p_nu_n(x_n)\to\iy}. \ee
By considering the maximum point of $u_n$, one can easily check $\SR\neq\emptyset$; see Section 5. Then Theorem \ref{thm-BM} tells us that after passing to a subsequence, $\SR\cap\Omega$ is an at most finite set. The problem is whether $\SR\cap\pa\Omega=\emptyset$ or not. When $W_n(x)\equiv1$, by the moving plane method one can prove that $p_nu_n$ is uniformly bounded in a small neighbourhood of $\pa\Omega$, and hence there is no boundary blow-up. However, due to the appearance of $W_n(x)$, the moving plane method is not applicable anymore. Here we use the induction method developed in \cite{LE-7} for the Lane-Emden equation in a small neighbourhood of $\pa\Omega$, and get $\SR\cap\pa\Omega=\emptyset$ by leading to a contradiction with $u_n|_{\pa\Omega}=0$. We point out that the induction method in \cite{LE-7} is inefficient at the places where $W_n(x)$ has zeros. Recall the zero set $\ZR$ of $W_n$ defined in \eqref{Z}. Once we obtain $\SR\subset\Omega$, we can apply Theorem \ref{thm-regular} near any point $a\in\SR\setminus\ZR$, and apply Theorem \ref{thm-singular} near any point $a\in\SR\cap\ZR$. Indeed, we obtain

\bt\lab{thm-1}
Let $u_n$ be a solution sequence of \eqref{equ-1-0}, and suppose $W_n(x)$ satisfies \eqref{con-1}-\eqref{con-2} with $\al_i>0$ for every $i=1,\cdots,m$. Then up to a subsequence, there exists a positive integer $k$ and different points $a_1,\cdots,a_k\in\Omega$ such that 
\begin{itemize}
\item[(i)]
The blow-up set $\SR$ of $p_nu_n$ is given by $\SR=\lbr{a_1,\cdots,a_k}$.
\item[(ii)] For small $r>0$, $\max_{B_r(a_i)}u_n\to \gamma_i\geq \sqrt{e}$ for all $1\leq i\leq k$. Furthermore, $\gamma_i=\sqrt{e}$ if $a_i\in\SR\setminus\ZR$.
\item[(iii)] For $t=0,1,2$, there holds
	\be\label{1234} p_n W_n(x)u_n(x)^{p_n-1+t}\to 8\pi e^{\frac{t}{2}} \Big(\sum_{a_i\in\SR\setminus\ZR}\delta_{a_i}+\sum_{a_j=q_{j'}\in\SR\cap\ZR}(1+\al_{j'})c_j^t\delta_{a_j} \Big),\ee
weakly in the sense of measures for some $c_{j}\ge1$.
\item[(iv)] For $a_j=q_{j'}\in\SR\cap \ZR$, $\gamma_j=\sqrt{e}$ and $c_j=1$ if $\al_{j'}\not\in\N$.
\end{itemize}\et

\br
\begin{itemize}
\item[(1)] For $\al_i\not\in\N$ for all $i=1,\cdots,m$, the existence of blow-up solutions of \eqref{equ-1-0} satisfying \eqref{1234} has been constructed by Esposito-Pistoia-Wei \cite{Reduction-1} via the finite-dimensional reduction method. In particular, their result shows that $\mathcal S\cap \mathcal Z\neq \emptyset$ happens for some solutions. Therefore, in general we can not expect $\mathcal S\cap \mathcal Z=\emptyset$ in Theorem \ref{thm-1}. In Theorem \ref{thm-1}, we prove in another direction that if $\al_i\not\in\N$ for all $i=1,\cdots,m$, then any solution sequence $u_n$ with bounded energy must behave the multi-point blow-up phenomena, and at any point the blow-up is simple. 

\item[(2)] For mean field equation with non-quantized singularity, i.e., $\al\not\in\N$, the profile of blow-up solutions has been given in Bartolucci-Tarantello \cite{SMF-2} and Bartolucci-Chen-Lin-Tarantello \cite{SMF-3}. They showed that the solution sequences develop multi-point blow-up and at each point 
the blow-up is simple. Our results are similar to theirs but different. As one can see we have no $\max u_n\to\iy$ but instead $\max u_n\to\sqrt e$; the energy of each bubble is dependent on the local maximum of $u_n$, which makes the analysis very different. 

\item[(3)] For the H\'enon equation
\be\lab{equ-henon} \begin{cases}
	-\Delta u_n=|x|^{2\al}u_n^{p_n},\quad u_n>0,\quad\text{in}~\Omega,\\
	u_n=0,\quad\text{on}~\pa \Omega,
\end{cases} \ee and some more general equations, the uniform bounded energy condition $\int_\Omega p_n|x|^{2\al}u_n^{p_n}\rd x\le C$ was proved to hold automatically in \cite{bound-1} for $\al>0$ and any simply connected domain $\Omega$ with $0\in\Omega$. 
\end{itemize}
\er

As an application of Theorem \ref{thm-1}, we study the ground states of the H\'enon equation \eqref{equ-henon} with $0\in\Omega$. 
Let $u_n(x)$ be a ground state (or called a least energy solution) of \eqref{equ-henon}, which by definition is a nontrivial solution of \eqref{equ-henon} such that the energy $\int_{\Omega}|x|^{2\alpha}|u_n|^{p_n+1}dx$ is smallest among all nontrivial solutions. It is standard to see that such ground state exists and is positive in $\Omega$ (up to a sign). We want to show that $0$ is not a blow-up point for the ground states.

When $\al=0$, the complete asymptotic behavior of the ground states as $p_n\to \infty$ was obtained in \cite{LE-5,LE-6,LE-8}, which says that the ground states behave as a single point blow-up (i.e. $k=1$ in \eqref{1-3}-\eqref{1-4}). 

For $\al>0$, Zhao \cite{Henon-1} proved some partial results for the ground state $u_n(x)$, which can be summarized as follows:
\begin{itemize}
\item For $\al>0$, 
\be\lab{tem-4-00}1\leq \liminf_{n\to\infty}\|u_n\|_{L^\infty(\Omega)}\leq \limsup_{n\to\infty}\|u_n\|_{L^\infty(\Omega)}\leq \sqrt{e},\ee
 \be\lab{tem-4}
	 \lim_{n\to\infty} p_n\int_{\Omega}|x|^{2\al}u_n^{p_n+1}\rd x= 8\pi e.
\ee
\item For $\al>e-1$, the ground state $u_n(x)$ behaves as at most two points blow-up, and $u_n(x)$ is not radially symmetric for $p_n$ large if $\Omega=B_r$ is an open ball.
\end{itemize} 
We want to improve these results and give a complete asymptotic behavior of the ground state $u_n(x)$ of the H\'enon equation \eqref{equ-henon} for any $\al>0$. 
To state our result, we introduce some notations. Recall the Green function $G(x,y)$
of $-\Delta$ in $\Omega$ with the Dirichlet boundary condition:
\be
	\left\{		\begin{aligned}
	&-\Delta_x G(x,y)=\delta_y 	\quad &\text{in}~\Omega,\\
	&G(x,y)=0					\quad &\text{on}~\pa\Omega,
	\end{aligned} 	\right.
\ee
It has the following form $$G(x,y)=-\f{1}{2\pi}\log|x-y|-H(x,y),\quad(x,y)\in\Omega\times\Omega,$$
where $H(x,y)$ is the regular part of $G(x,y)$.
It is well known that $H$ is a smooth function in $\Omega\times\Omega$, both
$G$ and $H$ are symmetric in $x$ and $y$. The Robin function of $\Omega$ is defined as
\be\lab{Robin} R(x):=H(x,x). \ee

\bt\lab{thm-0}
	Let $0\in\Omega$, $\al>0$ and $u_n$ be a ground state of the H\'enon equation \eqref{equ-henon}. Set $u_n(x_n)=\|u_n\|_{L^\infty(\Omega)}$. Then $u_n(x_n)\to\sqrt e$ and up to a subsequence, $x_n\to a\in\Omega\setminus\{0\}$, 
		$$p_nu_n\to 8\pi\sqrt e G(x,a),\quad\text{in}~\CR_{loc}^2(\overline{\Omega}\setminus\{a\}),$$
		\[ p_n|x|^{2\alpha}u_n(x)^{p_n-1+t}\to8\pi e^{\f{t}{2}}\delta_a, \quad t=0,1,2,\]
		$$\nabla \sbr{R(\cdot)-\f{1}{4\pi}\log|\cdot|^{2\al}}(a)=0,$$
	where $R(x)$ is the Robin function in \eqref{Robin}.
\et

\begin{remark} Theorem \ref{thm-0} improves those results in \cite{Henon-1}.
It is also interesting to compare Theorem \ref{thm-0} with some other results in the literature. Consider the H\'enon equation in general dimensions
\be \begin{cases}
	-\Delta u_n=|x|^{2\al}u_n^{p_n},\quad u_n>0,\quad\text{in}~\Omega,\\
	u_n=0,\quad\text{on}~\pa \Omega,
\end{cases} \ee
where $\Omega\subset\R^N$ is a smooth bounded domain. When $N\ge2$, the asymptotic behavior of ground states as $\al\to\iy$ was studied by Byeon-Wang in \cite{H-1,H-2}, where they proved that the ground states develop a boundary blow-up.
In another direction, when $N\ge3$, $\al>0$ is fixed, $\Omega=B_1$ and $p_n\to\frac{N+2}{N-2}$, Cao-Peng \cite{H-3} showed that the ground states also develop a boundary blow-up. However, Theorem \ref{thm-0} shows that there is no boundary blow-up for planar domains. Especially, when $\Omega$ is the unit ball, we know that the ground state of the H\'enon equation is not radially symmetric for $p_n$ large, since $x_n\to a\neq0$.
\end{remark}

The paper is organized as follows. In Section 2, we prove the Br\'ezis-Merle type result Theorem \ref{thm-BM}. In Sections 3 and 4, we study respectively the regular case and the singular case, and then prove Theorems \ref{thm-regular} and \ref{thm-singular}. In Section 5, we study the boundary value problem and prove Theorem \ref{thm-1}. Finally in Section 6, we study the ground states of the H\'enon equation. Throughout the paper, we denote by $C, C_0, C_1, \cdots$ to be positive constants independent of $n$ but may be different in different places.

\section{The Br\'ezis-Merle type result}

In this section, we follow Ren-Wei's idea \cite{LE-5,LE-6} to prove Theorem \ref{thm-BM}. Let $u_n$ be a solution sequence of \eqref{equ-BM} and denote
\be\lab{label-11}
	\bar u_n:=p_nu_n  \quad\text{and}\quad f_n:=p_nV_nu_n^{p_n}.
\ee
Then it holds
\be\lab{label-9}
	\begin{cases}
	-\Delta \bar u_n=f_n,\quad\text{in}~\Omega,\\
	\bar u_n>0,\quad\text{in}~\Omega,\\
	\bar u_n=0,\quad\text{on}~\pa\Omega.
	\end{cases}
\ee
Thanks to $\nm{f_n}_{L^1(\Omega)}\le C$, we may assume that 
	$$f_n\to\nu\quad\text{weakly in}~\MR(\Omega)~\text{as}~n\to\iy,$$
where $\MR(\Omega)$ is the space of Radon measures. Obviously $\nu(\Omega)<\iy$. 

For any $\delta>0$, we say a point $x_*\in\Omega$ to be \emph{a $\delta$-regular point} with respect to $\nu$, if there exists $\vp\in\CR_0(\Omega)$ satisfying $0\le\vp\le 1$, $\vp\equiv1$ near $x_*$ such that
	$$\int_\Omega\vp\rd\nu<\f{4\pi}{\tfrac{1}{e}+2\delta}.$$
Denote
	$$\Sigma_\nu(\delta):=\lbr{x\in\Omega:~x~\text{is not a $\delta$-regular point w.r.t. $\nu$}~}.$$
Before proceeding our discussion, we quote an $L^1$ estimates from \cite{MF-1}.
\bl[\cite{MF-1}]\lab{label-7}
	Let $u$ be a solution of
	$$\begin{cases}
		-\Delta u=f\quad\text{in}~\Omega,\\
		u=0\quad\text{on}~\pa\Omega,
	\end{cases}$$
	where $\Omega$ is a smooth bounded domain in $\R^2$. Then for any $0<\e<4\pi$, we have
		$$\int_{\Omega}\exp\sbr{\f{(4\pi-\e)|u(x)|}{\nm{f}_{L^1(\Omega)}}}\rd x\le \f{4\pi^2}{\e}(\mathrm{diam}~\Omega)^2.$$
\el

Now we give an equivalent characterization of the set $\Sigma_\nu(\delta)$.
\bl\lab{equiv-1}
	For any $\delta>0$ and $x_*\in\Omega$, we have that $x_*\in\Sigma_\nu(\delta)$ if and only if for any $R>0$ such that $B_R(x_*)\subset\Omega$, it holds $\nm{\bar u_n}_{L^\iy(B_R(x_*))}\to+\iy$ as $n\to\iy$. Consequently, $\Sigma_\nu(\delta)$ does not depend on the choice of $\delta$.
\el
\bp
First, take $x_*\not\in\Sigma_\nu(\delta)$, we want to prove that there exists $R_0>0$ such that $\nm{\bar u_n}_{L^\iy(B_{R_0}(x_*))}\le C$ as $n\to\iy$. Since $\nm{f_n}_{L^1(\Omega)}\le C$, by applying the elliptic $L^p$ estimate with the duality argument (cf. \cite{p-estimate}) to \eqref{label-9}, one gets that $\bar u_n$ are uniformly bounded in $W^{1,s}(\Omega)$ for any $1\le s<2$. In particular,
\be\lab{dd}
	\nm{\bar u_n}_{L^s(\Omega)}\leq C_s,\quad 1\leq s<2.
\ee
We claim that  there exist small $R_0>0$ and $\delta_0>0$ such that
\be\lab{claim-1}
	\nm{f_n}_{L^{1+\delta_0}(B_{2R_0}(x_*))}\le C,\quad\text{as}~n\to\iy.
\ee
Once \eqref{claim-1} is proved, we can apply the weak Harnack inequality (\cite[Theorem 8.17]{book-1}) to obtain
	$$\nm{\bar u_n}_{L^\iy(B_{R_0}(x_*))} \le C\sbr{\nm{\bar u_n}_{L^{3/2}(B_{2R_0}(x_*))}+\nm{f_n}_{L^{1+\delta_0}(B_{2R_0}(x_*))}}\le C. $$

Now we need to check the claim \eqref{claim-1}. Since $\f{\log x}{x}\le\f{1}{e}$ for any $x\in(0,+\iy)$, we obtain
	$$\log (p_n^{1/p_n}u_n(x))\le \f{1}{e}p_n^{1/p_n}u_n(x),\quad \forall x. $$
Therefore, for any $x\in\Omega$ and $\delta>0$
	$$f_n(x)=V_n(x)e^{p_n\log(p_n^{1/p_n}u_n(x))}\le Ce^{\f{1}{e}p_n^{1+1/p_n}u_n(x)}\le Ce^{(\f{1}{e}+\f{\delta}{2})\bar u_n(x)},\quad\text{for $n$ large}.$$
Since $x_*\not\in\Sigma_\nu(\delta)$, i.e. $x_*$ is a $\delta$-regular point, it follows from the definition of $\delta$-regular points that there exists $R_1>0$ such that $B_{2R_1}(x_*)\subset \Omega$ and
	$$\int_{B_{2R_1}(x_*)}f_n<\f{4\pi}{\tfrac{1}{e}+\delta}\quad\text{for $n$ large}.$$
Take $\bar u_n=\bar u_{n,1}+\bar u_{n,2}$ with $\bar u_{n,1}=0$ on the boundary $\pa B_{2R_1}(x_*)$ and $\bar u_{n,2}$ is harmonic in the ball $B_{2R_1}(x_*)$, i.e.
\be\label{harmonic-0}
\begin{cases}
-\Delta \bar u_{n,1}=f_n\quad\text{in }\; B_{2R_1}(x_*),\\
\bar u_{n,1}=0\quad\text{on }\;\pa B_{2R_1}(x_*),
\end{cases}
\qquad
\begin{cases}
-\Delta \bar u_{n,2}=0\quad\text{in }\; B_{2R_1}(x_*),\\
\bar u_{n,2}=\bar u_n\quad\text{on }\;\pa B_{2R_1}(x_*).
\end{cases}
\ee
 By the maximum principle, $\bar u_{n,1}>0$ and $\bar u_{n,2}>0$ in $B_{2R_1}(x_*)$. Applying Lemma \ref{label-7} to $\bar u_{n,1}$, we get
	$$\int_{B_{2R_1}(x_*)}\exp\sbr{\f{\ga \bar{u}_{n,1}(x)}{\nm{f_n}_{L^1(B_{2R_1}(x_*))}}}\rd x \le C_\ga,\quad\text{for any $\ga\in(0,4\pi)$}.$$
Note that $0<\bar u_{n,2}< \bar u_n$ in $B_{2R_1}(x_*)$. Then by the mean value theorem for harmonic functions and \eqref{dd}, we obtain
	$$\nm{\bar u_{n,2}}_{L^\iy(B_{R_1}(x_*))}\le C\nm{\bar u_{n,2}}_{L^1(B_{2R_1}(x_*))}\le C\nm{\bar u_n}_{L^1(B_{2R_1}(x_*))} \le C\nm{\bar u_{n}}_{L^1(\Omega)}\le C.$$
Take $\delta_0>0$ such that $\gamma:=4\pi(1+\delta_0)\frac{1+\frac{\delta}{2}e}{1+\delta e}<4\pi$. Then using the above estimates, we conclude that for $n$ large,
{\allowdisplaybreaks
\begin{align*}
		\int_{B_{R_1}(x_*)} f_n(x)^{1+\delta_0}\rd x
		&\le \int_{B_{R_1}(x_*)}C\exp\sbr{(1+\delta_0)(\tfrac{1}{e}+\tfrac{\delta}{2}) \bar u_n(x)} \rd x\\
		&\le C\int_{B_{R_1}(x_*)}\exp\sbr{(1+\delta_0)(\tfrac{1}{e}+\tfrac{\delta}{2}) \bar u_{n,1}(x)} \rd x\\
		&\le C\int_{B_{2R_1}(x_*)}\exp\sbr{(1+\delta_0)(\tfrac{1}{e}+\tfrac{\delta}{2}) \bar u_{n,1}(x)} \rd x\\
		&\le C\int_{B_{2R_1}(x_*)}\exp\sbr{4\pi(1+\delta_0)\f{1+\tfrac{\delta}{2}e}{1+\delta e}
				\f{\bar{u}_{n,1}(x)}{\nm{f_n}_{L^1(B_{2R_1}(x_*))}}}\rd x\\
		&=C\int_{B_{2R_1}(x_*)}\exp\sbr{
				\f{\gamma\bar{u}_{n,1}(x)}{\nm{f_n}_{L^1(B_{2R_1}(x_*))}}}\rd x\le C_{\gamma}.
	\end{align*}
	}%
Thus by choosing $R_0=R_1/2$, we finish the proof of the claim \eqref{claim-1}.

Finally, given any $x_*\in\Sigma_\nu(\delta)$, we claim that for any $R>0$, $\nm{\bar u_n}_{L^\iy(B_R(x_*))}\to+\iy$ as $n\to\iy$.
If not, then there exists $R_1>0$ such that up to a subsequence, $\nm{\bar u_n}_{L^\iy(B_{R_1}(x_*))}\le C$ as $n\to\iy$. Consequently,
	$$\int_{B_{R_1}(x_*)}f_n=\int_{B_{R_1}(x_*)}p_nV_n(x)u_n^{p_n}\le Cp_n\int_{B_{R_1}(x_*)}\sbr{\f{C}{p_n}}^{p_n} \to 0\quad\text{as $n\to\iy$.}$$
Thus by the definition of $\delta$-regular points, we obtain $x_*\not\in\Sigma_\nu(\delta)$, a contradiction. This finishes the proof.
\ep

\bc\lab{equiv-2}
	For any $\delta>0$, $\Sigma_\nu(\delta)\subset\Omega$ is an at most finite set.
\ec
\bp
Since $\nu(\lbr{x_*})\ge\f{4\pi}{\tfrac{1}{e}+2\delta}$ for every $x_*\in\Sigma_\nu(\delta)$,  it holds
	$$C\ge\nu(\Omega)\ge \f{4\pi}{\tfrac{1}{e}+2\delta}\#\Sigma_\nu(\delta),$$
which implies $\#\Sigma_\nu(\delta)<\iy$.
\ep

\bc\lab{label-10}
	For any compact subset $K\Subset\Omega\setminus\Sigma_\nu(\delta)$, it holds
		$$\nm{u_n}_{L^\iy(K)}\le \f{C_K}{p_n},\quad\text{for $n$ large}. $$
\ec
\bp
Given any compact subsets $K\Subset\Omega\setminus\Sigma_\nu(\delta)$, for any $x\in K$, we have $x\notin \Sigma_\nu(\delta)$, then it follows from Lemma \ref{equiv-1} that there exists $R_x>0$ such that 
	\[\nm{\bar u_n}_{L^\iy(B_{R_x}(x))}\le C_x,\quad\text{for $n$ large}.\]
From here and the finite covering theorem, we obtain
	$$\nm{\bar u_n}_{L^\iy(K)}\le C,\quad\text{for $n$ large}. $$
This implies
	$$\nm{u_n}_{L^\iy(K)}\le \f{C_K}{p_n},\quad\text{for $n$ large}. $$
Thus the proof is complete.
\ep

\bp[Proof of Theorem \ref{thm-BM}]
If $\Sigma_\nu(\delta)=\emptyset$, then Corollary \ref{label-10} implies that the alternative $(i)$ in Theorem \ref{thm-BM} holds. 

Thus we now suppose $\Sigma_\nu(\delta)\neq\emptyset$ and prove the alternative $(ii)$ in Theorem \ref{thm-BM} holds. Since $\Sigma_\nu(\delta)$ is a finite set and does not depend on the choice of $\delta$, we denote
\be\lab{label-12}
	\Sigma=\Sigma_\nu(\delta)=\lbr{a_1,\cdots,a_k}.
\ee
By Corollary \ref{label-10}, we know that $\nu=\sum_{i=1}^k\la_{a_i}\delta_{a_i}$. Since $\nu(\lbr{a_i})\ge\f{4\pi}{\tfrac{1}{e}+2\delta}$ for any $\delta>0$, we get $\la_{a_i}\ge4\pi e$, and hence $p_nV_n(x)u_n(x)^{p_n}\to\sum_{i=1}^k\la_{a_i}\delta_{a_i}$ weakly in the sense of measures in $\Omega$ with $\la_{a_i}\ge4\pi e$. Then by H\"older inequality, we get $p_nV_n(x)u_n(x)^{p_n-1}\to\sum_{i=1}^k\beta_{a_i}\delta_{a_i}$.

Choose $r_0>0$ such that
\be\lab{label-44} 
	B_{2r_0}(a_i)\subset\Omega\quad\text{and}\quad B_{2r_0}(a_i)\cap B_{2r_0}(a_j)=\emptyset,\quad\text{for}~i,j=1,\cdots,k,~i\neq j.
\ee
Define the local maximums $\ga_{n,i}$ and the local maximum points $x_{n,i}$ of $u_n$ by
\be\lab{notation-1}
	\ga_{n,i}=u_n(x_{n,i}):=\max_{B_{2r_0}(a_i)} u_n,\quad\text{for}~i=1,\cdots,k.
\ee

Recall the definition of $\bar u_n$ and $f_n$ in \eqref{label-11}, we have $-\Delta \bar u_n=f_n$ in $\Omega$. For any $i=1,\cdots,k$, it follows from \cite[Theorem 3.7]{book-1} that
\be  \max_{B_{2r_0}(a_i)}\bar u_n\le C(\max_{\pa B_{2r_0}(a_i)}\bar u_n+\max_{B_{2r_0}(a_i)}f_n), \ee
Since Corollary \ref{label-10} implies $\max_{\pa B_{2r_0}(a_i)}\bar u_n\le C$ and Lemma \ref{equiv-1} implies $\max_{B_{2r_0}(a_i)}\bar u_n\to+\iy$,  we have 
	$$\max_{B_{2r_0}(a_i)}f_n\to+\iy,$$
which, together with $\max_{B_{2r_0}(a_i)}f_n(x)\leq Cp_n \ga_{n,i}^{p_n}$, yields that up to a subsequence, $\ga_{n,i}\to\ga_i\ge1$. By Corollary \ref{label-10}, $u_n\to0$ in $L_{\loc}^\iy(\overline B_{2r_0}(a_i)\setminus\lbr{a_i})$, so $x_{n,i}\to a_i$ as $n\to\iy$. Finally, it is easy to see $\beta_{a_i}\ge\f{\la_{a_i}}{\ga_i}\ge\f{4\pi e}{\ga_i}$. This completes the proof.
\ep

\section{The regular case}
In this section, we prove Theorem \ref{thm-regular}. Let $u_n$ be a solution sequence of \eqref{equ-regular}. Without loss of generality, we may assume the radius $r=1$. Suppose $0$ is the only blow-up point of $p_nu_n$ in $B_1$ and $V_n(x)$ satisfies \eqref{con-regular2}.

Let $x_n$ be a maximum point of $u_n$ in $B_1$, i.e.
	$$u_n(x_n):=\max_{\overline B_1}u_n,$$
then \eqref{con-regular1} implies $p_nu_n(x_n)\to\iy$ and $x_n\to0$. Define the scaling parameter $\mu_n>0$ by
\be \mu_n^{-2}:=p_nV_n(x_n)u_n(x_n)^{p_n-1}, \ee
and the scaling function by
	\be\label{55}v_n(x):=p_n\sbr{\f{u_n(x_n+\mu_nx)}{u_n(x_n)}-1} \quad\text{for}~x\in D_n:=\f{B_1-x_n}{\mu_n}.\ee
It is easy to see that $v_n$ satisfies
\be\lab{label-1}
	\begin{cases}
	-\Delta v_n=\f{V_n(x_n+\mu_nx)}{V_n(x_n)}\sbr{1+\f{v_n}{p_n}}^{p_n}\quad\text{in}~D_n,\\
	v_n(0)=0=\max_{D_n}v_n,
	\end{cases}
\ee
and
\be\label{5-0}0<1+\f{v_n(x)}{p_n}=\f{u_n(x_n+\mu_nx)}{u_n(x_n)}\leq 1\quad \text{in }D_n.\ee

\bl\lab{converge-regular}
	After passing to a subsequence, it hold $\mu_n\to0$, $u_n(x_n)\to\ga\in[1,\iy)$ and $v_n\to U_0(x)=-2\log(1+\f{1}{8}|x|^2)$ in $\CR_{loc}^2(\R^2)$.
\el
\bp
Suppose $\mu_n\not\to0$, then up to subsequence we may assume $u_n(x_n)\le\sbr{\f{C}{p_n}}^{\f{1}{p_n-1}}$ for some constant $C>0$. Thus it holds $0\le -\Delta (p_nu_n)\le C$, which together with $\displaystyle\max_{\pa B_1}p_nu_n\le C$ implies $\displaystyle\max_{B_1}p_nu_n\le C$. This is a contradiction with that $0$ is a blow-up point of $p_nu_n$. So $\mu_n\to0$ and hence $u_n(x_n)\to\ga\ge 1$. 

Now we prove $\ga<\iy$. Recall \cite[Proposition 2.7]{LE-4} that there is $C>0$ independent of $x\in\Omega$ and $p$ such that
	 $$\nm{G(x,\cdot)}_{L^{p}(\Omega)}^{p}\leq Cp^{p+1},\quad \text{for $p>1$ large}.$$
Then by the Green's representation formula and H\"older inequality,
{\allowdisplaybreaks
\begin{align*}
	u_n(x_n)
	&=\int_{B_1} G(x_n,y)V_n(y)u_n(y)^{p_n}\rd y-\int_{\pa B_1}\f{\pa G(x_n,y)}{\pa\nu}u_n(y)\rd s_y \\
	&\le C\nm{G(x_n,\cdot)}_{L^{2p_n+1}(B_1)}\sbr{\int_{B_1}V_n(y)^{1+\f{1}{2p_n}}u_n(y)^{p_n+\f{1}{2}}}^{\f{2p_n}{2p_n+1}}\\
	&\qquad+\f{C}{p_n}\int_{\pa B_1}\abs{\f{\pa G(x_n,y)}{\pa\nu}}\rd s_y\\
	&\le C(2p_n+1)^{\f{2p_n+2}{2p_n+1}}u_n(x_n)^{\f{p_n}{2p_n+1}}\sbr{\int_{B_1}V_n(y)u_n(y)^{p_n}}^{\f{2p_n}{2p_n+1}}+\f{C}{p_n}\\
	&\le C(2p_n+1)^{\f{2p_n+2}{2p_n+1}}u_n(x_n)^{\f{p_n}{2p_n+1}}\sbr{\f{C}{p_n}}^{\f{2p_n}{2p_n+1}}+\f{C}{p_n}\\
	&\le C\sbr{u_n(x_n)^{\f{p_n}{2p_n+1}}+\f{1}{p_n}},\quad\text{for $p_n>1$ large enough},
\end{align*}
}%
so $\limsup_{n\to\iy} u_n(x_n)\le C$, i.e., $\ga<\iy$.

For any $R>0$, $B_R\subset D_n$ for $n$ large. Like \eqref{harmonic-0} we let
	$$v_n=\vp_n+\psi_n\quad\text{in}~B_R,$$
with $-\Delta \vp_n=-\Delta v_n$ in $B_R$ and $\psi_n=v_n$ on $\pa B_R$. Thanks to \eqref{label-1}-\eqref{5-0}, we see that $|-\Delta v_n|\le C$ in $D_n$ for some constant $C>0$. Then by the standard elliptic theory, we obtain that $\vp_n$ is uniformly bounded in $B_R$. Since $\psi_n=v_n-\vp_n$, we know that $\psi_n$ is harmonic in $B_R$ and bounded from above. By the Harnack inequality, we see that if $\inf_{B_R}\psi_n\to-\iy$, then $\sup_{B_R}\psi_n\to-\iy$ as $n\to\iy$, which contradicts with $\psi_n(0)=-\vp_n(0)\ge-C$. So $\psi_n$ and hence $v_n$ is uniformly bounded in $B_R$. After passing to a subsequence, the standard elliptic theory implies that 
	$$v_n\to U_0 \quad\text{in}~\CR_{loc}^2(\R^2)~\text{as}~n\to\iy,$$
and \eqref{label-1} implies
\be\lab{label-2}
\begin{cases}
	-\Delta U_0=e^{U_0}\quad\text{in}~\R^2,\\
	U_0(0)=0=\max_{\R^2}U_0.
		\end{cases}
\ee
Moreover, by Fatou's Lemma,
	\begin{align*}\int_{\R^2}e^{U_0}\rd x&\le \liminf_{n\to\iy}\int_{D_n}\f{V_n(x_n+\mu_nx)}{V_n(x_n)}\Big(1+\f{v_n}{p_n}\Big)^{p_n}\rd x\\
	&=\liminf_{n\to\iy}\f{p_n}{u_n(x_n)}\int_{B_1}V_n(x)u_n(x)^{p_n}\rd x\le C.\end{align*}
Since $U_0(0)=0$, by the classification result due to Chen and Li \cite{classification-1} we obtain
	$$U_0(x)=-2\log \sbr{1+\f{1}{8}|x|^2},$$
and $\int_{\R^2}e^{U_0}\rd x=8\pi$. 
\ep

We introduce the local Pohozaev identity.
\bl
	Suppose $u$ satisfies
	$$\begin{cases} -\Delta u=V(x)u^p,\quad\text{in}~\Omega,\\ u>0,\quad\text{in}~\Omega,\end{cases}$$
	then for any $y\in\R^2$ and any subset $\Omega'\subset\Omega$, it holds
	\be\lab{pho-1} \begin{aligned} 
		&\quad \f{1}{p+1}\int_{\Omega'}\sbr{2V(x)+\abr{\nabla V(x),x-y}}u(x)^{p+1}\rd x \\
		&=\int_{\pa\Omega'}\abr{\nabla u(x),\nu(x)}\abr{\nabla u(x),x-y}-\f{1}{2}|\nabla u(x)|^2\abr{x-y,\nu(x)}\rd s_x \\
		&\quad +\f{1}{p+1}\int_{\pa\Omega'} V(x)u(x)^{p+1}\abr{x-y,\nu(x)}\rd s_x, 
	\end{aligned} \ee
where $\nu(x)$ denotes the outer normal vector of $\partial \Omega'$ at $x$.
\el
\bp
By direct computations, we have
	$$-\Delta u(x)\cdot\abr{\nabla u(x),x-y}=-\operatorname{div}\sbr{\nabla u(x)\abr{\nabla u(x),x-y}-\f{1}{2}|\nabla u(x)|^2 (x-y)},$$
and 
	$$\begin{aligned} 
		V(x)u(x)^p\cdot\abr{\nabla u(x),x-y}&=\f{1}{p+1}\operatorname{div}\sbr{V(x)u(x)^{p+1}(x-y)}\\
		&\quad -\f{1}{p+1}(2V(x)+\abr{\nabla V(x),x-y})u(x)^{p+1}.
	\end{aligned}$$
Then multiplying $-\Delta u=V(x)u^p$ with $\abr{\nabla u(x),x-y}$, integrating on $\Omega'$ and using the divergence theorem, we obtain \eqref{pho-1}.
\ep

By \eqref{con-regular1} we have that for any compact subset $K\Subset\overline B_1\setminus\{0\}$,
\be\lab{label-36} \nm{p_nu_n}_{L^\iy(K)}\le C_K. \ee

\bl\lab{converge-regular2}
	It holds
		$$p_nu_n(x)\to 8\pi\ga G_1(x,0)+\psi(x),\quad\text{in}~\CR_{loc}^2(\overline B_1\setminus\{0\})~\text{as}~n\to\iy. $$
	where $\gamma$ is given in Lemma \ref{converge-regular}, $\psi\in\CR^2(\overline B_1)$ is a harmonic function, and $G_1(x,y)$ denotes the Green function of $-\Delta$ in $B_1$ with the Dirichlet boundary condition.
\el
\bp
Like \eqref{harmonic-0} we set $u_n=\phi_n+\psi_n$ with $\phi_n=0$ on $\pa B_1$ and $\psi_n$ is harmonic in $B_1$. Since $\psi_n=u_n=O(\f{1}{p_n})$ on $\pa B_1$, it follows from the standard elliptic theory that up to a subsequence, $p_n\psi_n\to\psi$ in $\CR^2(\overline B_1)$. Since $p_n\psi_n$ is harmonic, so is $\psi$.

Take $d\in(0,1)$ and any compact subset $K\Subset\overline B_1\setminus\{0\}$. Applying the Green's representation formula to $\phi_n$ and using \eqref{con-regular1}-\eqref{con-regular2}, we get that for any $x\in K$,
$$\begin{aligned}
	p_n\phi_n(x)
	&=\int_{B_1} G_1(x,y)p_nV_n(y)u_n(y)^{p_n} \rd y\\
	&=\int_{B_d} G_1(x,y)p_nV_n(y)u_n(y)^{p_n} \rd y+o_n(1)\int_{B_1\setminus B_d} G_1(x,y)\rd y\\
	&\to\sigma_0 G_1(x,0),\quad\text{uniformly for $x\in K$ as}~n\to\iy,
\end{aligned}$$
where 
\be\lab{label-43}   \sigma_0:=\lim_{d\to0}\lim_{n\to\iy}\int_{B_d}p_nV_n(x)u_n(x)^{p_n}\rd x. \ee
Again by the Green's representation formula, a similar argument implies
	$$\nabla_x(p_n\phi_n)(x)=\int_\Omega \nabla_x G_1(x,y)p_nV_n(y)u_n(y)^{p_n} \rd y\to \sigma_0 \nabla_xG_1(x,0).$$
Thus $p_nu_n(x)\to \sigma_0 G_1(x,0)+\psi(x)$ in $\CR_{loc}^1(\overline B_1\setminus\{0\})$. From here and $-\Delta (p_nu_n)=p_nV_n(x)u_n^{p_n}\to 0$ in $L_{\loc}^\iy(\overline B_1\setminus\{0\})$ and $p_nu_n\to\psi$ in $\CR^2(\pa B_1)$, it follows from the standard elliptic estimates that
$$p_nu_n(x)\to \sigma_0 G_1(x,0)+\psi(x),\quad\text{in}~\CR_{loc}^2(\overline B_1\setminus\{0\})~\text{as}~n\to\iy. $$

It remains to prove $\sigma_0=8\pi\ga$. Since
$$\begin{aligned}
	\int_{B_d}p_nV_n(x)u_n(x)^{p_n}\rd x
	&=u_n(x_n)\int_{\frac{B_{d}-x_n}{\mu_n}}\f{V_n(x_n+\mu_nx)}{V_n(x_n)}\sbr{1+\f{v_n(x)}{p_n}}^{p_n}\rd x\\
	&\ge \ga\int_{\R^2}e^{U_0}\rd x+o_n(1)=8\pi\ga+o_n(1),
\end{aligned}$$
we get $\sigma_0\ge 8\pi\ga$. On the other hand,
applying the Pohozaev identity \eqref{pho-1} with $y=0$, $\Omega'=B_d$, $V=V_n$ and $u=u_n$, we obtain
\be\lab{label-46} \begin{aligned}
	&\quad \f{p_n^2}{p_n+1}\int_{B_d}\sbr{2V_n(x)+\abr{\nabla V_n(x),x}}u_n(x)^{p_n+1}\rd x \\
	&=\int_{\pa B_d}\abr{p_n\nabla u_n(x),\nu(x)}\abr{p_n\nabla u_n(x),x}-\f{1}{2}|p_n\nabla u_n(x)|^2\abr{x,\nu(x)}\rd s_x \\
	&\quad +\f{p_n^2}{p_n+1}\int_{\pa B_d} V_n(x)u_n(x)^{p_n+1}\abr{x,\nu(x)}\rd s_x. 
\end{aligned} \ee
Note that for $x\in B_1\setminus\{0\}$, we have
\be\lab{label-45}  p_n\nabla u_n(x)\to\sigma_0\nabla_x G(x,0)+\nabla\psi(x)=-\f{\sigma_0}{2\pi}\f{x}{|x|^2}+O(1).  \ee
Using \eqref{label-36} and \eqref{label-45}, we obtain (note $\nu(x)=\f{x}{|x|}$ on $\pa B_d$)
	$$\lim_{n\to\infty}\text{RHS of \eqref{label-46}}=\f{\sigma_0^2}{4\pi}+O(d).$$
	From here and \eqref{label-46}, we conclude
	\be\label{label-48}
	\lim_{d\to 0}\lim_{n\to\infty}p_n\int_{B_d}\sbr{2V_n(x)+\abr{\nabla V_n(x),x}}u_n(x)^{p_n+1}\rd x=\f{\sigma_0^2}{4\pi}.
	\ee
	Since $V_n$ satisfies \eqref{con-regular2}, we have
	$$\abs{ p_n\int_{B_d}\abr{\nabla V_n(x),x}u_n(x)^{p_n+1}\rd x }\le Cd\int_{B_d}p_nV_n(x)u_n(x)^{p_n}\rd x\le C d,$$
which together with \eqref{label-48} and \eqref{label-43} implies
\[\f{\sigma_0^2}{8\pi}=\lim_{d\to0}\lim_{n\to\iy} p_n\int_{B_d}V_n(x)u_n(x)^{p_n+1}\rd x\leq \lim_{n\to\infty}u_n(x_n) \sigma_0=\gamma\sigma_0,\]
so $\sigma_0\le 8\pi\ga$. This proves $\sigma_0=8\pi\ga$.
\ep

For the scaling function $v_n$ defined in \eqref{55}, we need the following decay estimates, which will be used to apply the Dominated Convergence Theorem.
\bl\lab{decay-regular1} 
For any $\eta\in(0,4)$, there exist small $r_\eta>0$, large $R_\eta>1$, $n_\eta>1$ and constant $C_\eta>0$ such that
\be  v_n(x)\le \eta\log\f{1}{|x|}+C_\eta\quad \text{and} \quad \abs{v_n(x)}\le C_\eta(1+\log|x|),  \ee
for any $2R_\eta\le |x|\le \f{r_\eta}{\mu_n}$ and $n\ge n_\ga$.
\el
\bp
By Lemma \ref{converge-regular} we have 
	$$v_n(x)\to U_0(x)=-2\log\sbr{1+\f{1}{8}|x|^2}\quad \text{in}~\CR_{loc}^2(\R^2).$$
Moreover, Lemma \ref{converge-regular2} tells $\sigma_0= 8\pi\ga$ with $\sigma_0$ defined by \eqref{label-43}. 

Applying the Green's representation formula, we have for any $x\in D_n$,
$$\begin{aligned}
	u_n(x_n+\mu_nx)
	&=\int_{B_1} G_1(x_n+\mu_nx,y)V_n(y)u_n(y)^{p_n}\rd y-\int_{\pa B_1} \f{\pa G_1(x_n+\mu_nx,y)}{\pa\nu_y}u_n(y)\rd s_y\\
	&=\f{u_n(x_n)}{p_n}\int_{D_n} G_1(x_n+\mu_nx,x_n+\mu_nz) \f{V_n(x_n+\mu_nz)}{V_n(x_n)} \sbr{1+\f{v_n(z)}{p_n}}^{p_n}\rd z\\
	&\quad -\int_{\pa B_1} \f{\pa G_1(x_n+\mu_nx,y)}{\pa\nu_y}u_n(y)\rd s_y.
\end{aligned}$$
Then it follows from \eqref{55} that
$$\begin{aligned}
	v_n(x)&=-p_n+ \int_{D_n} G_1(x_n+\mu_nx,x_n+\mu_nz)\f{V_n(x_n+\mu_nz)}{V_n(x_n)}\sbr{1+\f{v_n(z)}{p_n}}^{p_n}\rd z\\
	&\quad -\f{p_n}{u_n(x_n)}\int_{\pa B_1} \f{\pa G_1(x_n+\mu_nx,y)}{\pa\nu_y}u_n(y)\rd s_y.
\end{aligned}$$
Since $v_n(0)=0$ and $G_1(z,y)=-\f{1}{2\pi}\log|z-y|-H_1(z,y)$, we have
{\allowdisplaybreaks
\begin{align*}
	&\quad v_n(x)\\
	&=v_n(x)-v_n(0)\\
	&=\int_{D_n} \mbr{G_1(x_n+\mu_nx,x_n+\mu_nz)-G_1(x_n,x_n+\mu_nz)}\f{V_n(x_n+\mu_nz)}{V_n(x_n)}\sbr{1+\f{v_n(z)}{p_n}}^{p_n}\rd z\\
		&\quad -\f{p_n}{u_n(x_n)}\int_{\pa B_1} \sbr{\f{\pa G_1(x_n+\mu_nx,y)}{\pa\nu_y}-\f{\pa G_1(x_n,y)}{\pa\nu_y}}u_n(y)\rd s_y\\
	&=\f{1}{2\pi}\int_{D_n} \log\f{|z|}{|z-x|}\f{V_n(x_n+\mu_nz)}{V_n(x_n)}\sbr{1+\f{v_n(z)}{p_n}}^{p_n}\rd z\\
		&\quad -\int_{D_n} \mbr{H_1(x_n+\mu_nx,x_n+\mu_nz)-H_1(x_n,x_n+\mu_nz)}\f{V_n(x_n+\mu_nz)}{V_n(x_n)}\sbr{1+\f{v_n(z)}{p_n}}^{p_n}\rd z\\
		&\quad -\f{p_n}{u_n(x_n)}\int_{\pa B_1} \sbr{\f{\pa G_1(x_n+\mu_nx,y)}{\pa\nu_y}-\f{\pa G_1(x_n,y)}{\pa\nu_y}}u_n(y)\rd s_y\\
	&=:\Rmnum{1}(x)+\Rmnum{2}(x)+\Rmnum{3}(x).
\end{align*}
}%
Since $H_1(x,y)$ is smooth in $B_1\times B_1$ and $\nabla G_1(x,y)$ is bounded for $|x-y|\ge c>0$, we have that for $|x|\le \f{r_\eta}{\mu_n}$ with small $r_\eta<\f{1}{2}$ to be chosen later, 
$$\begin{aligned}
	\Rmnum{2}(x)
	&=O(1)\int_{D_n} \f{V_n(x_n+\mu_nz)}{V_n(x_n)}\sbr{1+\f{v_n(z)}{p_n}}^{p_n}\rd z\\
	&=O(1)\f{p_n}{u_n(x_n)}\int_{B_1} V_n(y)u_n(y)^{p_n}\rd y=O(1),
\end{aligned}$$
$$\begin{aligned}
	\Rmnum{3}(x)=O(1)\int_{\pa B_1} \abs{\f{\pa G_1(x_n+\mu_nx,y)}{\pa\nu_y}}+\abs{\f{\pa G_1(x_n,y)}{\pa\nu_y}}\rd y=O(1).
\end{aligned}$$

For any fixed $\eta\in(0,4)$, let $\e=\f{2\pi}{3}(4-\eta)>0$ and take $R_\eta>1$ large such that $\int_{B_{R_\eta}(0)}e^{U_0}>\int_{\R^2}e^{U_0}-\f{\e}{2}=8\pi-\f{\e}{2}$, where $U_0(z)=-2\log\sbr{1+\f{1}{8}|z|^2}$. Then from $v_n\to U_0$ we get that for $n$ large,
\be\lab{label-54}   \int_{B_{R_\eta}} \f{V_n(x_n+\mu_nz)}{V_n(x_n)}\sbr{1+\f{v_n(z)}{p_n}}^{p_n}\rd z \ge \int_{B_{R_\eta}}e^{U_0(z)}\rd z-\f{\e}{2}> 8\pi-\e. \ee
From $\sigma_0=8\pi\ga$, we see that
	\begin{align*}&\lim_{r\to0}\lim_{n\to\iy}\int_{\lbr{|z|\le\f{2r}{\mu_n}}}\f{V_n(x_n+\mu_nz)}{V_n(x_n)}\sbr{1+\f{v_n(z)}{p_n}}^{p_n}\rd z\\
	=&\lim_{r\to0}\lim_{n\to\iy}\frac{p_n}{u_n(x_n)}\int_{B_{2r}}V_n(y)u_n(y)^{p_n}dy=\frac{\sigma_0}{\gamma}= 8\pi.\end{align*}
Thus we can choose $r_\eta\in (0,\frac12)$ small such that for $n$ large,
\be\lab{label-55}   \int_{\lbr{|z|\le\f{2r_\eta}{\mu_n}}}\f{V_n(x_n+\mu_nz)}{V_n(x_n)}\sbr{1+\f{v_n(z)}{p_n}}^{p_n}\rd z\le 8\pi+\e.\ee
By \eqref{label-54}-\eqref{label-55} we obtain
\begin{equation}\label{65}
\int_{\lbr{R_{\eta}\leq |z|\le\f{2r_\eta}{\mu_n}}}\f{V_n(x_n+\mu_nz)}{V_n(x_n)}\sbr{1+\f{v_n(z)}{p_n}}^{p_n}\rd z\le 2\e.
\end{equation}

Fix any $2R_\eta\le |x|\le \f{r_\eta}{\mu_n}$, to compute the integral $\Rmnum{1}(x)$, we divide $ D_n$ into four domains $ D_n=\sum_{i=1}^4D_{n,i}$ and divide the integral $\Rmnum{1}(x)$ into four terms $\Rmnum{1}(x)=\sum_{i=1}^4\Rmnum{1}_{D_{n,i}}(x)$, where $D_{n,1}=\lbr{z\in D_n:~|z|\le R_\eta}$, $D_{n,2}=\lbr{z\in D_n:~|z|\ge \f{2r_\eta}{\mu_n}}$ and
	$$D_{n,3}=\lbr{z\in D_n:~R_\eta\le|z|\le\f{2r_\eta}{\mu_n},~ |z|\le2|z-x|\le3|z|},$$
	$$D_{n,4}=\lbr{z\in D_n:~R_\eta\le|z|\le\f{2r_\eta}{\mu_n},~|z|\ge2|z-x|~\text{or}~2|z-x|\ge 3|z|}.$$

If $z\in D_{n,1}$, then $|x|\ge 2|z|$, $|x-z|\geq \frac12 |x|$ and hence
	$$\log\f{|z|}{|z-x|}\le \log\f{2R_\ga}{|x|}\le 0.$$
From here and \eqref{label-54}, we have
\begin{align}\lab{label-60} 
\Rmnum{1}_{D_{n,1}}(x)&\le \f{1}{2\pi}\log\f{2R_\eta}{|x|}\int_{D_{n,1}} \f{V_n(x_n+\mu_nz)}{V_n(x_n)}\sbr{1+\f{v_n(z)}{p_n}}^{p_n}\rd z \\
&\le (4-\f{\e}{2\pi})\log\f{1}{|x|}+C.\nonumber\end{align}
On the other hand, since  $v_n(x)\le 0$ and
	$$0\ge\log\f{1}{|z-x|}\ge\log\f{2}{3|x|},\quad\text{for }z\in D_{n,1},$$
we get 
\begin{align}\lab{label-61} 
	\Rmnum{1}_{D_{n,1}}(x)
	&\ge \f{1}{2\pi}\log\f{2}{3|x|}\int_{D_{n,1}} \f{V_n(x_n+\mu_nz)}{V_n(x_n)}\sbr{1+\f{v_n(z)}{p_n}}^{p_n}\rd z -C\int_{D_{n,1}}\abs{ \log|z|}\rd z\nonumber\\
	&\ge (4+\frac{\varepsilon}{2\pi})\log\f{1}{|x|}-C.
\end{align}
Note that $|z|\geq 2|x|$ for $z\in D_{n,2}$. Then it is easy to see that
	$$\log\f{2}{3}\le\log\f{|z|}{|z-x|}\le\log2,\quad\text{for }\;z\in D_{n,2}\cup D_{n,3},$$
which implies
\begin{align}\lab{label-62} 
	\abs{\Rmnum{1}_{D_{n,2}}(x)+\Rmnum{1}_{D_{n,3}}(x)}&\le C\int_{ D_n} \f{V_n(x_n+\mu_nz)}{V_n(x_n)}\sbr{1+\f{v_n(z)}{p_n}}^{p_n}\rd z\\
	&= \frac{Cp_n}{u_n(x_n)}\int_{\Omega}V_n(y)u_n(y)^{p_n}\rd y\le C.\nonumber\end{align}
Finally for $z\in D_{n,4}$, it holds $2\le |z|\le 2|x|$.  Then we see from \eqref{65} that
\be\lab{label-56} \begin{aligned}
	0&\le \f{1}{2\pi}\int_{D_{n,4}} \log|z|\f{V_n(x_n+\mu_nz)}{V_n(x_n)}\sbr{1+\f{v_n(z)}{p_n}}^{p_n}\rd z\\
	&\le \f{1}{2\pi}\log(2|x|)\int_{D_{n,4}}\f{V_n(x_n+\mu_nz)}{V_n(x_n)}\sbr{1+\f{v_n(z)}{p_n}}^{p_n}\rd z\\
	&\le \f{\e}{\pi}\log|x|+C,
\end{aligned}\ee
Furthermore, by $v_n(x)\le 0$, we get
\be\lab{label-57} \begin{aligned}
	0&\le \f{1}{2\pi}\int_{\lbr{z\in D_{n,4}:~|z-x|\le1}} \log\f{1}{|z-x|}\f{V_n(x_n+\mu_nz)}{V_n(x_n)}\sbr{1+\f{v_n(z)}{p_n}}^{p_n}\rd z\\
	&\le C\int_{\lbr{|z-x|\le1}}\log\f{1}{|z-x|}\rd x\le C.
\end{aligned}\ee
While for $z\in\lbr{z\in D_{n,4}:~|z-x|\ge1}$, it holds
	$$\log\f{1}{3|x|}\le \log\f{1}{|z-x|}\le 0,$$
and hence
\be\lab{label-59} \begin{aligned}
	0\ge&\f{1}{2\pi}\int_{\lbr{z\in D_{n,4}:~|z-x|\ge1}} \log\f{1}{|z-x|}\f{V_n(x_n+\mu_nz)}{V_n(x_n)}\sbr{1+\f{v_n(z)}{p_n}}^{p_n}\rd z\\
	&\ge \f{1}{2\pi}\log\f{1}{3|x|}\int_{D_{n,4}}\f{V_n(x_n+\mu_nz)}{V_n(x_n)}\sbr{1+\f{v_n(z)}{p_n}}^{p_n}\rd z\\
	&\ge \f{\e}{\pi}\log\f{1}{3|x|}.
\end{aligned}\ee
Combining \eqref{label-56}-\eqref{label-59}, we get 
\be\lab{label-63} \f{\e}{\pi}\log\f{1}{3|x|}\le \Rmnum{1}_{D_{n,4}}(x)\le \f{\e}{\pi}\log|x|+C. \ee

By \eqref{label-60},\eqref{label-61},\eqref{label-62},\eqref{label-63} and $\Rmnum{2}(x)=O(1)$, $\Rmnum{3}(x)=O(1)$, we finally get
	$$\abs{v_n(x)}\le C(1+\log|x|),$$
and
	$$v_n(x)\le \mbr{4-\f{3\e}{2\pi}}\log\f{1}{|x|}+C=\eta\log\f{1}{|x|}+C,$$
for any $2R_\eta\le |x|\le \f{r_\eta}{\mu_n}$ and some constant $C>0$. This completes the proof.
\ep
\br\lab{decay-regular2}
For any $\eta\in(0,4)$, Lemma \ref{decay-regular1} implies
	$$\sbr{1+\f{v_n(x)}{p_n}}^{p_n}=e^{p_n\log(1+\f{v_n(x)}{p_n})}\le e^{v_n(x)} \le \f{C_\eta}{|x|^\eta},\quad\text{\it for }\;2R_\eta\le |x|\le \f{r_\eta}{\mu_n}.$$
Meanwhile, since $v_n\to U_0$ in $\CR_{loc}^2(\R^2)$, we have $\sbr{1+\f{v_n(x)}{p_n}}^{p_n}\le C$ for $|x|\le 2R_\eta$ and $n$ large. 
Therefore,
\be    0\le \sbr{1+\f{v_n(x)}{p_n}}^{p_n}\le \f{C_\eta}{1+|x|^\eta}, \quad\forall|x|\leq \f{r_\eta}{\mu_n}. \ee
Similarly, we have
\be 	\abs{v_n(x)}\le C_\eta\log\sbr{2+|x|}, \quad\forall|x|\leq \f{r_\eta}{\mu_n}. \ee
\er

As a direct application of the above decay estimates, we have
\bl\lab{tem-2}
	It holds $\ga=\sqrt e$.
\el
\bp Take $\eta=3$ in Lemma \ref{decay-regular1} and Remark \ref{decay-regular2} and let $r=r_3/2$. Let $G_{r}(x,y)$ denote the Green function of $-\Delta$ in $B_{r}$ with the Dirichlet boundary condition.
By the Green's representation formula, we have
$$\begin{aligned}
	& u_n(x_n)=\int_{B_{r}} G_{r}(x_n,y)V_n(y)u_n(y)^{p_n}\rd y-\int_{\pa B_{r}}\f{\pa G_{r}(x_n,y)}{\pa\nu}u_n(y)\rd s_y \\
	&=\f{u_n(x_n)}{p_n}\int_{\frac{B_{r}-x_n}{\mu_n}} G_{r}(x_n,x_n+\mu_ny)\f{V_n(x_n+\mu_ny)}{V_n(x_n)}\sbr{1+\f{v_n(y)}{p_n}}^{p_n}\rd y+O(\f{1}{p_n}),\\
\end{aligned}$$
so
	\[\f{1}{p_n}\int_{\frac{B_{r}-x_n}{\mu_n}} G_{r}(x_n,x_n+\mu_ny)\f{V_n(x_n+\mu_ny)}{V_n(x_n)}\sbr{1+\f{v_n(y)}{p_n}}^{p_n}\rd y=1+O(\f{1}{p_n}).\]
On the other hand, by Remark \ref{decay-regular2}, for any $y\in \frac{B_{r}-x_n}{\mu_n}$, we have $|y|\le \f{r_3}{\mu_n}$ for $n$ large and so
\be\lab{33-25} 		0\le\sbr{1+\f{v_n(y)}{p_n}}^{p_n}\le \f{C}{1+|y|^3}.\ee
Then by applying the Dominated Convergence Theorem, we get
	\[ \lim_{n\to\infty}\int_{\frac{B_{r}-x_n}{\mu_n}} \f{V_n(x_n+\mu_ny)}{V_n(x_n)}\sbr{1+\f{v_n(y)}{p_n}}^{p_n}\rd y=\int_{\R^2}e^{U_0}\rd x=8\pi, \]
\begin{align*}
	&\lim_{n\to\infty}\int_{\frac{B_{r}-x_n}{\mu_n}} \sbr{\f{1}{2\pi}\log|y|+H_{r}(x_n,x_n+\mu_ny)}\f{V_n(x_n+\mu_ny)}{V_n(x_n)}\sbr{1+\f{v_n(y)}{p_n}}^{p_n}\rd y\\
	&\qquad=\int_{\R^2} \sbr{\f{1}{2\pi}\log|y|+H_{r}(0,0)} e^{U_0(y)}\rd y=C <\infty.
\end{align*}
From here, $G_{r}(x,y)=-\f{1}{2\pi}\log|x-y|-H_{r}(x,y)$ and $\mu_n^{-2}=p_nV_n(x_n)u_n(x_n)^{p_n-1}$, we have
$$\begin{aligned}
&\quad 1+O(\f{1}{p_n})\\
	&= \f{1}{p_n}\int_{\frac{B_{r}-x_n}{\mu_n}} G_{r}(x_n,x_n+\mu_ny)\f{V_n(x_n+\mu_ny)}{V_n(x_n)}\sbr{1+\f{v_n(y)}{p_n}}^{p_n}\rd y\\
	&=-\f{1}{2\pi}\f{\log\mu_n}{p_n}\int_{\frac{B_{r}-x_n}{\mu_n}} \f{V_n(x_n+\mu_ny)}{V_n(x_n)}\sbr{1+\f{v_n(y)}{p_n}}^{p_n}\rd y+O(\frac1{p_n})\\
	&=\f{1}{4\pi}\sbr{\f{\log p_n+\log V_n(x_n)}{p_n}+\f{p_n-1}{p_n}\log u_n(x_n)}(8\pi+o_n(1))+o_n(1)\\
	&=2\log \ga+o_n(1).
\end{aligned}$$
Thus $2\log \ga=1$, i.e., $\ga=\sqrt e$.
\ep

Now we are ready to prove Theorem \ref{thm-regular}.
\bp[Proof of Theorem \ref{thm-regular}]
It has been proved in Lemma \ref{tem-2} that $\displaystyle\max_{\overline B_1}u_n\to \sqrt e$. Since $0$ is the only blow-up point of $p_nu_n$ in $B_1$, we see that 
	$$p_nV_n(x)u_n(x)^{p_n-1+k}\to \beta_k\delta_0, \quad\text{for }k=0,1,2, $$
	weakly in the sense of measures.
For any small $r>0$, it follows from the Dominated Convergence Theorem that
{\allowdisplaybreaks
\begin{align*}
	&\int_{B_r}p_nV_n(x)u_n(x)^{p_n-1+k}\rd x\\
	=&u_n(x_n)^{k}\int_{\frac{B_{r}-x_n}{\mu_n}}\f{V_n(x_n+\mu_nx)}{V_n(x_n)}\sbr{1+\f{v_n(x)}{p_n}}^{p_n-1+k}\rd x
	\to \gamma^{k}\int_{\R^2}e^{U_0}\rd x=8\pi e^{\frac{k}{2}},
\end{align*}
}%
Thus $\beta_k=8\pi e^{\frac{k}{2}}$ for $k=0,1,2$. This completes the proof.
\ep

\section{The singular case}

In this section, we prove Theorem \ref{thm-singular}. Let $u_n$ be a solution sequence of \eqref{equ-singular}. Without loss of generality, we may assume the radius $r=1$. Suppose $0$ is the only blow-up point of $p_nu_n$ in $B_1$ and $V_n(x)$ satisfies \eqref{con-singular2}.

Let $x_n$ be a maximum point of $u_n$ in $B_1$, i.e.
\be\lab{max-1}   u_n(x_n)=\max_{\overline B_1}u_n,  \ee
then \eqref{con-singular1} implies $p_nu_n(x_n)\to\iy$ and $x_n\to0$. We claim that
\be\label{4-42}    \mu_n^{-2-2\al}:=p_nV_n(x_n)u_n(x_n)^{p_n-1}\to\iy. \ee
Indeed, if $p_nV_n(x_n)u_n(x_n)^{p_n-1}\not\to\iy$, then up to a subsequence, we have $u_n(x_n)\le\sbr{\f{C}{p_n}}^{\f{1}{p_n-1}}$ for some constant $C>0$. Thus it holds $0\le -\Delta (p_nu_n)\le C$, which together with $\displaystyle\max_{\pa B_1}p_nu_n\le C$ implies $\displaystyle\max_{B_1}p_nu_n\le C$. This is a contradiction with that $0$ is a blow-up point of $p_nu_n$. So $p_nV_n(x_n)u_n(x_n)^{p_n-1}\to\iy$ and hence $\liminf_{n\to\iy}u_n(x_n)\ge 1$. Following the approach in Lemma \ref{converge-regular}, we may assume that 
\be\lab{max}   u_n(x_n)=\max_{\overline B_1}u_n\to\ga\in[1,\iy). \ee
Up to a subsequence, we denote
\begin{align}\label{tem-60-beta}
\beta_k:=\lim_{n\to\infty} \int_{B_1}p_n|x|^{2\al}V_n(x)u_n(x)^{p_n-1+k}\rd x,\quad k=0,1,2.
\end{align}
Then
\be\beta_2\leq \gamma \beta_1\leq \gamma^2\beta_0.\ee
Furthermore, for any $0<d<1$, since \eqref{con-singular1} gives $\sup_{B_1\setminus B_{d}}p_nu_n\le C$, we have
\[ \lim_{n\to\infty} \int_{B_1\setminus B_d}p_n|x|^{2\al}V_n(x)u_n(x)^{p_n-1+k}\rd x=0,\]
so
\begin{align}\label{tem-60-beta1}
\beta_k=\lim_{d\to 0}\lim_{n\to\infty} \int_{B_d}p_n|x|^{2\al}V_n(x)u_n(x)^{p_n-1+k}\rd x,\quad k=0,1,2.
\end{align}

\bl\lab{tem-203} We have 
	\be\label{beta12}\beta_1^2=8\pi(1+\alpha)\beta_2,\quad \beta_1\leq 8\pi(1+\alpha)\gamma.\ee
\el

\bp Let $G_1(x,y)$ denotes the Green function of $-\Delta$ in $B_1$ with the Dirichlet boundary condition.
Exactly as in Lemma \ref{converge-regular2}, we get
\be p_nu_n(x)\to \beta_1 G_1(x,0)+\psi(x),\quad\text{in}~\CR_{loc}^2(\overline B_1\setminus\{0\})~\text{as}~n\to\iy, \ee
where $\psi\in\CR^2(\overline B_1)$ is a harmonic function.
Consequently, for $x\in B_1\setminus\{0\}$, we have
\be\lab{tem-22}  p_n\nabla u_n(x)\to\beta_1\nabla_x G(x,0)+\nabla\psi(x)=-\f{\beta_1}{2\pi}\f{x}{|x|^2}+O(1).  \ee
Applying the Pohozaev identity \eqref{pho-1} with $y=0$, $\Omega'=B_d$, $V=|x|^{2\al}V_n(x)$ and $u=u_n$, and by using \eqref{con-singular1} and \eqref{tem-22}, we obtain
\be\lab{tem-23} \begin{aligned}
	&\quad \lim_{n\to\infty}\f{p_n^2}{p_n+1}\int_{B_d}|x|^{2\al}\mbr{(2+2\al)V_n(x)+\abr{\nabla V_n(x),x}}u_n(x)^{p_n+1}\rd x \\
	&=\f{\beta_1^2}{4\pi}+O(d).
\end{aligned} \ee
Since $V_n$ satisfies \eqref{con-singular2}, we have
	$$\abs{ p_n\int_{B_d}|x|^{2\al}\abr{\nabla V_n(x),x}u_n(x)^{p_n+1}\rd x }\le Cd\int_{B_d}p_n|x|^{2\al}V_n(x)u_n(x)^{p_n}\rd x\le C d.$$
From here and \eqref{tem-23}, we deduce 
	$$\f{\beta_1^2}{8\pi(1+\al)}=\lim_{d\to0}\lim_{n\to\iy}\int_{B_d}p_n|x|^{2\al}V_n(x)u_n(x)^{p_n+1}\rd x=\beta_2\le \ga\beta_1,$$
namely \eqref{beta12} holds.
\ep

Since $|x_n|\to0$ and $p_nV_n(x_n)u_n(x_n)^{p_n-1}\to\iy$, we need compare their convergence rates to analyse the values of $\beta_k$.

\subsection{A special case.} \
In this section, we assume that
\be\lab{assume-1}  p_n|x_n|^{2+2\al}V_n(x_n)u_n(x_n)^{p_n-1}\le C. \ee
Define the scaling function
	$$v_n(x):=p_n\sbr{\f{u_n(x_n+\mu_nx)}{u_n(x_n)}-1} \quad\text{for}~x\in D_n:=\f{B_1-x_n}{\mu_n}.$$
It is easy to see that $v_n$ satisfies
\be\lab{tem-20}
	\begin{cases}
	-\Delta v_n=\abs{x+\f{x_n}{\mu_n}}^{2\al}\f{V_n(x_n+\mu_nx)}{V_n(x_n)}\sbr{1+\f{v_n}{p_n}}^{p_n}\quad\text{in}~D_n,\\
	v_n(0)=0=\max_{D_n}v_n.
	\end{cases}
\ee
Since \eqref{assume-1} implies $|\frac{x_n}{\mu_n}|\leq C$, up to a subsequence we have $\f{x_n}{\mu_n}\to x_\iy$ for some $x_\iy\in\R^2$. Then by following the approach of Lemma \ref{converge-regular}, we obtain $v_n\to U_\al$ in $\CR_{loc}^2(\R^2)$, where $U_\al$ satisfies
\be\begin{cases}
	-\Delta U_\al=|x+x_\iy|^{2\al}e^{U_\al}\quad\text{in}~\R^2,\\
	U_\al(0)=0=\max_{\R^2}U_\al,\\
	\int_{\R^2}e^{U_{\al}}\rd x\le C.
\end{cases}\ee
By the classification result due to Prajapat and Tarantello \cite{classification-2}, we obtain 
\be\lab{tem-24}   U_\al(z)=-2\log \sbr{1+\f{1}{8(1+\al)^2}|(z+z_\iy)^{1+\al}-z_\iy^{1+\al}|^2},\quad z\in\C,\ee
where $z_\iy\in\mathbb{C}$ is the complex notation of $x_\iy$. Moreover, \[\int_{\R^2}|x+x_\iy|^{2\al}e^{U_\al}\rd x=8\pi(1+\al). \]

\bl\lab{converge-singular2}
	Suppose \eqref{assume-1} holds, then $\beta_1=8\pi(1+\al)\ga$.
\el

\bp
By Fatou's Lemma, for any $d\in (0,1)$, 
$$\begin{aligned}
	&\quad \lim_{n\to\iy}\int_{B_d}p_n|x|^{2\al}V_n(x)u_n(x)^{p_n}\rd x\\
	&=u_n(x_n)\int_{\frac{B_{d}-x_n}{\mu_n}}\abs{x+\f{x_n}{\mu_n}}^{2\al}\f{V_n(x_n+\mu_nx)}{V_n(x_n)}\sbr{1+\f{v_n(x)}{p_n}}^{p_n}\rd x\\
	&\ge \ga\int_{\R^2}|x+x_\iy|^{2\al}e^{U_\al}\rd x+o_n(1)=8\pi(1+\al)\ga+o_n(1),
\end{aligned}$$
we get $\beta_1\ge 8\pi(1+\al)\ga$. Together with \eqref{beta12}, 
we obtain $\beta_1= 8\pi(1+\al)\ga$.
\ep

Since the following lemma is similar to Lemma \ref{decay-regular1}, we sketch the proof and only emphasize the different places.

\bl\lab{decay-singular1} 
Suppose \eqref{assume-1} holds, then for any $\eta\in(0,4(1+\al))$, there exist small $r_\eta>0$, large $R_\eta>1$, $n_\eta>1$ and constant $C_\eta>0$ such that
\be  v_n(x)\le \eta\log\f{1}{|x|}+C_\eta\quad \text{and} \quad \abs{v_n(x)}\le C_\eta(1+\log|x|),  \ee
for any $2R_\eta\le |x|\le \f{r_\eta}{\mu_n}$ and $n\ge n_\ga$.
\el
\bp
Exactly as in Lemma \ref{decay-regular1}, we have 
\be\lab{tem-39}\begin{aligned}
	v_n(x)
	&=\f{1}{2\pi}\int_{D_n} \log\f{|z|}{|z-x|}\abs{z+\f{x_n}{\mu_n}}^{2\al}\f{V_n(x_n+\mu_nz)}{V_n(x_n)}\sbr{1+\f{v_n(z)}{p_n}}^{p_n}\rd z+O(1)\\
	&=:\Rmnum{1}(x)+O(1)
\end{aligned}\ee
for $|x|\le \f{r_\eta}{\mu_n}$ with small $r_\eta<\f{1}{2}$ to be chosen later.

For any fixed $\eta\in(0,4(1+\al))$, let $\e=\f{2\pi}{3+4\al}[4(1+\al)-\eta]>0$ and take $R_\eta>0$ such that 
	$$\int_{B_{R_\eta}}|z+x_\iy|^{2\al}e^{U_\al(z)}\rd z>\int_{\R^2}|z+x_\iy|^{2\al}e^{U_\al(z)}\rd z-\f{\e}{2}=8\pi(1+\al)-\f{\e}{2},$$
where $U_\al$ is given in \eqref{tem-24}. Then from $v_n\to U_\al$ we get that for $n$ large,
\begin{align}\lab{tem-26}  & \int_{B_{R_\eta}} \abs{z+\f{x_n}{\mu_n}}^{2\al}\f{V_n(x_n+\mu_nz)}{V_n(x_n)}\sbr{1+\f{v_n(z)}{p_n}}^{p_n}\rd z \\
\ge &\int_{B_{R_\eta}}|z+x_\iy|^{2\al}e^{U_\al}\rd z-\f{\e}{2}> 8\pi(1+\al)-\e. \nonumber\end{align}
From $\beta_1=8\pi(1+\al)\ga$ with $\beta_1$ satisfying \eqref{tem-60-beta1}, we see that
	$$\lim_{r\to0}\lim_{n\to\iy}\int_{\lbr{|z|\le\f{2r}{\mu_n}}}\abs{z+\f{x_n}{\mu_n}}^{2\al}\f{V_n(x_n+\mu_nz)}{V_n(x_n)}\sbr{1+\f{v_n(z)}{p_n}}^{p_n}\rd z= 8\pi(1+\al).$$
Thus we can choose $r_\eta\in (0,\frac12)$ small such that for $n$ large,
\be\lab{tem-27}   \int_{\lbr{|z|\le\f{2r_\eta}{\mu_n}}}\abs{z+\f{x_n}{\mu_n}}^{2\al}\f{V_n(x_n+\mu_nz)}{V_n(x_n)}\sbr{1+\f{v_n(z)}{p_n}}^{p_n}\rd z\le 8\pi(1+\al)+\e,\ee
and consequently,
\be\lab{tem-27-0}   \int_{\lbr{R_{\eta}\leq |z|\le\f{2r_\eta}{\mu_n}}}\abs{z+\f{x_n}{\mu_n}}^{2\al}\f{V_n(x_n+\mu_nz)}{V_n(x_n)}\sbr{1+\f{v_n(z)}{p_n}}^{p_n}\rd z\le 2\e.\ee

Fix any $2R_\eta\le |x|\le \f{r_\eta}{\mu_n}$, to compute the integral $\Rmnum{1}(x)$, we divide $ D_n$ into the same four domains $ D_n=\sum_{i=1}^4D_{n,i}$ as in Lemma \ref{decay-regular1} and divide the integral $\Rmnum{1}(x)$ into four terms $\Rmnum{1}(x)=\sum_{i=1}^4\Rmnum{1}_{D_{n,i}}(x)$.
Then as in Lemma \ref{decay-regular1}, we obtain
\be\lab{tem-31}  [4(1+\al)+\f{\e}{2\pi}]\log\f{1}{|x|}-C\le \Rmnum{1}_{D_{n,1}}(x)\le [4(1+\al)-\f{\e}{2\pi}]\log\f{1}{|x|}+C,\ee
\be\lab{tem-32} \abs{\Rmnum{1}_{D_{n,2}}(x)+\Rmnum{1}_{D_{n,3}}(x)}\le C.\ee
Finally for $z\in D_{n,4}$, it holds $2\le |z|\le 2|x|$.  Then it follows from \eqref{tem-27-0} that
\be\lab{tem-33}  \begin{aligned}
	0&\le \f{1}{2\pi}\int_{D_{n,4}} \log|z|\abs{z+\f{x_n}{\mu_n}}^{2\al}\f{V_n(x_n+\mu_nz)}{V_n(x_n)}\sbr{1+\f{v_n(z)}{p_n}}^{p_n}\rd z\\
	&\le \f{1}{2\pi}\log(2|x|)\int_{D_{n,4}}\abs{z+\f{x_n}{\mu_n}}^{2\al}\f{V_n(x_n+\mu_nz)}{V_n(x_n)}\sbr{1+\f{v_n(z)}{p_n}}^{p_n}\rd z\\
	&\le \f{\e}{\pi}\log|x|+C.
\end{aligned}\ee
Furthermore, by $v_n(x)\le 0$, we get
$$ \begin{aligned}
	0&\le \f{1}{2\pi}\int_{\lbr{z\in D_{n,4}:~|z-x|\le\f{1}{|x|^{2\al}}}} \log\f{1}{|z-x|}\abs{z+\f{x_n}{\mu_n}}^{2\al}\f{V_n(x_n+\mu_nz)}{V_n(x_n)}\sbr{1+\f{v_n(z)}{p_n}}^{p_n}\rd z\\
	&\le C(1+|x|^{2\al})\int_{\lbr{|z-x|\le\f{1}{|x|^{2\al}}}}\log\f{1}{|z-x|}\rd z\\
	&\le C\f{(1+|x|^{2\al})(1+\log|x|)}{|x|^{4\al}}\le C.
\end{aligned}$$
While for $z\in\lbr{z\in D_{n,4}:~|z-x|\ge\f{1}{|x|^{2\al}}}$, it holds
	$$\log\f{1}{3|x|}\le \log\f{1}{|z-x|}\le 2\al\log|x|,$$
and hence
$$\begin{aligned}
	&~\f{1}{2\pi}\int_{\lbr{z\in D_{n,4}:~|z-x|\ge\f{1}{|x|^{2\al}}}} \log\f{1}{|z-x|}\abs{z+\f{x_n}{\mu_n}}^{2\al}\f{V_n(x_n+\mu_nz)}{V_n(x_n)}\sbr{1+\f{v_n(z)}{p_n}}^{p_n}\rd z\\
	&\le \f{\al}{\pi}\log|x| \int_{D_4}\abs{z+\f{x_n}{\mu_n}}^{2\al}\f{V_n(x_n+\mu_nz)}{V_n(x_n)}\sbr{1+\f{v_n(z)}{p_n}}^{p_n}\rd z \\
	&\le \f{2\al\e}{\pi}\log|x|,
\end{aligned}$$
\begin{align}\lab{tem-30-2} 
	\f{1}{2\pi}\int_{\lbr{z\in D_{n,4}:~|z-x|\ge\f{1}{|x|^{2\al}}}}& \log\f{1}{|z-x|}\abs{z+\f{x_n}{\mu_n}}^{2\al}\f{V_n(x_n+\mu_nz)}{V_n(x_n)}\sbr{1+\f{v_n(z)}{p_n}}^{p_n}\rd z\nonumber\\
	\ge &\f{1}{2\pi}\log \f{1}{3|x|}\int_{D_{n,4}}\abs{z+\f{x_n}{\mu_n}}^{2\al}\f{V_n(x_n+\mu_nz)}{V_n(x_n)}\sbr{1+\f{v_n(z)}{p_n}}^{p_n}\rd z\\
	\ge &\f{\e}{\pi}\log \f{1}{3|x|}.\nonumber
\end{align}
Combining \eqref{tem-33}-\eqref{tem-30-2}, we get 
\be\lab{tem-40} \f{\e}{\pi}\log\f{1}{3|x|}\le \Rmnum{1}_{D_{n,4}}(x)\le \f{(1+2\al)\e}{\pi}\log|x|+C. \ee

Finally, from \eqref{tem-39},\eqref{tem-31},\eqref{tem-32} and \eqref{tem-40}, we finally get
	$$\abs{v_n(x)}\le C(1+\log|x|),$$
and
	$$v_n(x)\le \mbr{4(1+\al)-\f{(3+4\al)\e}{2\pi}}\log\f{1}{|x|}+C=\eta\log\f{1}{|x|}+C,$$
for any $2R_\eta\le |x|\le \f{r_\eta}{\mu_n}$ and some constant $C>0$. This completes the proof.
\ep
\br\lab{decay-singular2}
Similarly as Remark \ref{decay-regular2}, we have that for any $\eta\in(0,4(1+\al))$, there exists $C_\eta>0$ such that
\be    0\le \sbr{1+\f{v_n(x)}{p_n}}^{p_n}\le \f{C_\eta}{1+|x|^\eta}, \quad\forall|x|\leq \f{r_{\eta}}{\mu_n}, \ee
\be 	\abs{v_n(x)}\le C_\eta\log\sbr{2+|x|}, \quad\forall|x|\leq \f{r_{\eta}}{\mu_n}. \ee
\er

\bl\lab{tem-41}
	Suppose \eqref{assume-1} holds, then $\ga=\sqrt e$ and 
		$$p_n|x|^{2\al}V_n(x)u_n(x)^{p_n-1+k}\to 8\pi(1+\al)e^{\frac{k}{2}}\delta_0,\quad\text{for }k=0,1,2 $$
	weakly in the sense of measures.
\el
\bp  Take $\eta=3(1+\alpha)$ in Lemma \ref{decay-singular1} and Remark \ref{decay-singular2}, and let $r=r_{3(1+\alpha)}/2$.
Then the same argument as Lemma \ref{tem-2} implies
	\[\f{1}{p_n}\int_{\frac{B_{r}-x_n}{\mu_n}} G_{r}(x_n,x_n+\mu_ny)\abs{y+\f{x_n}{\mu_n}}^{2\al}\f{V_n(x_n+\mu_ny)}{V_n(x_n)}\sbr{1+\f{v_n(y)}{p_n}}^{p_n}\rd y=1+O(\frac1{p_n}).\]
Since $|y|\leq \frac{r_{3(1+\alpha)}}{\mu_n}$ for $y\in \frac{B_{r}-x_n}{\mu_n}$, we see from Remark \ref{decay-singular2} that
\be		0\le\sbr{1+\f{v_n(y)}{p_n}}^{p_n}\le \f{C}{1+|y|^{3(1+\al)}},\quad\text{\it for }\;y\in \frac{B_{r}-x_n}{\mu_n},\ee
so it follows from the Dominated Convergence Theorem that
	\begin{align*} &\lim_{n\to\infty}\int_{\frac{B_{r}-x_n}{\mu_n}} \abs{y+\f{x_n}{\mu_n}}^{2\al}\f{V_n(x_n+\mu_ny)}{V_n(x_n)}\sbr{1+\f{v_n(y)}{p_n}}^{p_n}\rd y\\
	=&\int_{\R^2}|y+x_\iy|^{2\al}e^{U_\al}\rd y=8\pi(1+\al), \end{align*}
\begin{align*}
	&\lim_{n\to\infty}\int_{\frac{B_{r}-x_n}{\mu_n}} \sbr{\f{1}{2\pi}\log|y|+H_{r}(x_n,x_n+\mu_ny)}\abs{y+\f{x_n}{\mu_n}}^{2\al}\f{V_n(x_n+\mu_ny)}{V_n(x_n)}\sbr{1+\f{v_n(y)}{p_n}}^{p_n}\rd y\\
	&\qquad=\int_{\R^2} \sbr{\f{1}{2\pi}\log|y|+H_{r}(0,0)} |y+x_\iy|^{2\al}e^{U_\al(y)}\rd y=C <\infty.
\end{align*}
From here, $G_{r}(x,y)=-\f{1}{2\pi}\log|x-y|-H_{r}(x,y)$ and $\mu_n^{-2-2\al}=p_nV_n(x_n)u_n(x_n)^{p_n-1}$, we have
{\allowdisplaybreaks
\begin{align*}
&\quad 1+O(\frac1{p_n})\\
	&=\f{1}{p_n}\int_{\frac{B_{r}-x_n}{\mu_n}} G_{r}(x_n,x_n+\mu_ny)\abs{y+\f{x_n}{\mu_n}}^{2\al}\f{V_n(x_n+\mu_ny)}{V_n(x_n)}\sbr{1+\f{v_n(y)}{p_n}}^{p_n}\rd y\\
	&=-\f{1}{2\pi}\f{\log\mu_n}{p_n}\int_{\frac{B_{r}-x_n}{\mu_n}} \abs{y+\f{x_n}{\mu_n}}^{2\al}\f{V_n(x_n+\mu_ny)}{V_n(x_n)}\sbr{1+\f{v_n(y)}{p_n}}^{p_n}\rd y+O(\f{1}{p_n})\\
	&=\frac{1}{4\pi(1+\alpha)}\sbr{\f{\log p_n+\log V_n(x_n)}{p_n}+\f{p_n-1}{p_n}\log u_n(x_n)}(8\pi(1+\alpha)+o_n(1))+o_n(1)\\
	&=2\log \ga+o_n(1).
\end{align*}
}%
Thus $2\log \ga=1$, i.e., $\ga=\sqrt e$.

Since $0$ is the only blow-up point of $p_nu_n$ in $B_1$, we see that 
	$$p_nV_n(x)u_n(x)^{p_n-1+k}\to \beta_k\delta_0, \quad \text{for }k=0,1,2, $$
	weakly in the sense of measures, where $\beta_k$ are given by \eqref{tem-60-beta}-\eqref{tem-60-beta1}.
For any small $r>0$, again by the Dominated Convergence Theorem, we get
$$\begin{aligned}
	&\quad\int_{B_r}p_nV_n(x)|x|^{2\al}u_n(x)^{p_n-1+k}\rd x\\
	&=u_n(x_n)^{k}\int_{\frac{B_{r}-x_n}{\mu_n}}\abs{x+\f{x_n}{\mu_n}}^{2\al}\f{V_n(x_n+\mu_nx)}{V_n(x_n)}\sbr{1+\f{v_n(x)}{p_n}}^{p_n-1+k}\rd x\\
	&\to \ga^{k}\int_{\R^2}|x+x_\iy|^{2\al}e^{U_\al}\rd x=8\pi (1+\al)e^{\frac{k}{2}}.
\end{aligned}$$
Therefore, $\beta_k=8\pi(1+\al)e^{\frac{k}{2}}$, and the proof is complete.
\ep

\subsection{The general case.} In this section, we do not assume 
the estimate \eqref{assume-1}, so the previous arguments in Section 4.1 do not work and different ideas are needed. 
We begin with the following decomposition result, whose proof is inspired by \cite[Proposition 1.4]{SMF-1}.

\bpr\lab{decomposition}
Let $u_n$ satisfy the assumptions of Theorem \ref{thm-singular}.  Then along a subsequence, one of the following alternatives holds.
\begin{itemize}
\item[(i)] 	Either there exists $\e_0\in(0,\f{1}{2})$ such that
	\be\lab{bound} 	\sup_{B_{2\e_0}}p_n|x|^{2+2\al}u_n(x)^{p_n-1}\le C, \ee
\item[(ii)]	 or there exist $\e_0'\in(0,\f{1}{2})$ and $l\geq 1$ sequences $\{z_{n,i}\}_{n\ge1}\subset B_1\setminus\{0\}$, $i=1,\cdots,l$, such that
	\be\lab{tem-15}  \lim_{n\to\infty} z_{n,i}=0,\quad \liminf_{n\to\infty}|z_{n,i}|^{\f{2+2\al}{p_n-1}}u_n(z_{n,i})\geq 1, \ee
		\be\lab{tem-15-0}   \lim_{n\to\infty} p_n|z_{n,i}|^{2+2\al}u_n(z_{n,i})^{p_n-1}=\iy, \ee
	\be\lab{tem-16}   p_n|x|^{2+2\al}u_n(x)^{p_n-1}\le C,\;\text{for}~x\in\lbr{y\in B_1:~|y|\le 2\e_0'|z_{n,1}|\;\text{or}\; |y|\ge \f{1}{2\e_0'}|z_{n,l}|}, \ee
	and in case $l\ge 2$, then $\f{|z_{n,i}|}{|z_{n,i+1}|}\to0$ as $n\to\iy$ and
	\be\lab{tem-17}   p_n|x|^{2+2\al}u_n(x)^{p_n-1}\le C,\;\text{for}~x\in\bigcup_{i=1}^{l-1}\lbr{y\in B_1:~\f{1}{2\e_0'}|z_{n,i}|\le |y|\le 2\e_0'|z_{n,i+1}|}. \ee
\end{itemize}
\epr
\bp
We devide the proof into several steps.

\vskip0.1in
{\bf Step 1.} Assume there exists $\{z_n\}\subset B_1$ such that $p_n|z_n|^{2+2\al}u_n(z_n)^{p_n-1}\to\iy$, we prove that
$\lim_{n\to\infty}z_n=0$, $\liminf_{n\to\infty}|z_{n}|^{\f{2+2\al}{p_n-1}}u_n(z_{n})\geq 1$ and 
\be\lab{tem-3}   \limsup_{n\to\iy}\int_{B_{\delta|z_n|}(z_n)}p_n|x|^{2\al}V_n(x)u_n(x)^{p_n}\rd x\ge 4\pi e,\quad\text{for every}~\delta>0. \ee

Indeed, $p_n|z_n|^{2+2\al}u_n(z_n)^{p_n-1}\to\iy$ implies \[\liminf_{n\to\iy}u_n(z_n)\ge1,\quad \liminf_{n\to\infty}|z_{n}|^{\f{2+2\al}{p_n-1}}u_n(z_{n})\geq 1,\] so $p_nu_n(z_n)\to \infty$. Since  \eqref{con-singular1} says that $0$ is the only blow-up point of $p_nu_n$ in $B_1$, we obtain $z_n\to 0$. In particular, this argument shows that once \eqref{tem-15-0} holds, then \eqref{tem-15} holds.

To prove \eqref{tem-3}, we let
\be\lab{al_n}  v_n(x):=|z_n|^{\al_n}u_n(|z_n|x), \quad\text{with}~\al_n:=\f{2+2\al}{p_n-1}. \ee
Since $|z_n|\to0$  and  $p_n|z_n|^{2+2\al}u_n(z_n)^{p_n-1}\to\iy$, we get
	$$1\ge|z_n|^{\al_n}\ge p_n^{-\f{1}{p_n-1}}\f{1}{u_n(z_n)}\ge\f{1}{\ga}+o(1),$$
where $\ga$ is given in \eqref{max}. Then
\be   \int_{B_{\f{1}{|z_n|}}}p_n|x|^{2\al}V_n(|z_n|x)v_n(x)^{p_n}\rd x=|z_n|^{\al_n} \int_{B_1}p_n|x|^{2\al}V_n(x)u_n(x)^{p_n}\rd x\le C. \ee
From \eqref{equ-singular}, we see that 
\be\label{4-32}\begin{cases}
	-\Delta v_n=|x|^{2\al}V_n(|z_n|x)v_n^{p_n},\quad\text{in}~B_{\f{1}{|z_n|}},\\
	p_nv_n(\f{z_n}{|z_n|})\to\iy.
\end{cases}\ee
Take a subsequence so that $\f{z_n}{|z_n|}$ converges to some point $x_0$ in the unit circle. Then $p_nv_n$ admits a blow-up point at $x_0$ and around it the function $W_n(x)=|x|^{2\al}V_n(|z_n|x)$ is uniformly bounded from above. Then by using Theorem \ref{thm-BM} for $v_n$ in any open bounded domain $\Omega$ containing $x_0$, we see that for any small $\delta>0$,
	$$ \limsup_{n\to\iy}\int_{B_{\delta}(\f{z_n}{|z_n|})}p_n|x|^{2\al}V_n(|z_n|x)v_n(x)^{p_n}\rd x\ge 4\pi e.$$ 
	A simple change of variables leads to 
\begin{align*}	
	&\limsup_{n\to\iy}\int_{B_{\delta |z_{n}|}(z_n)}p_n|x|^{2\al}V_n(x)u_n(x)^{p_n}\rd x\\
	=&\limsup_{n\to\iy}|z_n|^{-\alpha_n}\int_{B_{\delta}(\f{z_n}{|z_n|})}p_n|x|^{2\al}V_n(|z_n|x)v_n(x)^{p_n}\rd x\geq 4\pi e,
\end{align*}
namely \eqref{tem-3} holds. This proves Step 1.

\vskip0.1in
{\bf Step 2.} Suppose the alternative (i) does not hold for every $\e_0\in(0,\f{1}{2})$, we prove that there exist $\e_0'\in(0,\f{1}{2})$ and a sequence $\{z_{n,1}\}$ such that
\be\lab{tem-5}  \lim_{n\to\infty} z_{n,1}=0,\quad  \lim_{n\to\infty} p_n|z_{n,1}|^{2+2\al}u_n(z_{n,1})^{p_n-1}=\iy, \ee
and 
\be\lab{tem-6}   p_n|x|^{2+2\al}u_n(x)^{p_n-1}\le C,\quad\text{for}~x\in\lbr{y\in B_1:~|y|\le 2\e_0'|z_{n,1}|}. \ee

Indeed, since \eqref{bound} does not hold for every $\e_0\in(0,\f{1}{2})$,  up to a subsequence there is $z_n\in {B}_{1}$ such that
	$$p_n|z_n|^{2+2\al}u_n(z_n)^{p_n-1}=\sup_{B_{1/2}}p_n|x|^{2+2\al}u_n(x)^{p_n-1}\to\iy.$$
Then by Step 1 we have $z_n\to0$ and
	$$\limsup_{n\to\iy}\int_{B_{\delta|z_n|}(z_n)}p_n|x|^{2\al}V_n(x)u_n(x)^{p_n}\rd x\ge 4\pi e,\quad\text{for every}~\delta>0.$$
Let $v_n$ be defined by \eqref{al_n}, then $v_n$ satisfies \eqref{4-32}. There are the same alternatives for $v_n$. If there exists some $\e_0'\in(0,\f{1}{2})$ such that
\be\lab{tem-7} 	\sup_{B_{2\e_0'}}p_n|x|^{2+2\al}v_n(x)^{p_n-1}\le C, \ee
then 
\[\sup_{B_{2\e_0'|z_n|}}p_n|x|^{2+2\al}u_n(x)^{p_n-1}=\sup_{B_{2\e_0'}}p_n|x|^{2+2\al}v_n(x)^{p_n-1}\le C,\]
so setting $z_{n,1}=z_n$ we get \eqref{tem-5}-\eqref{tem-6} and we are done. Otherwise, for any $r\in (0,\frac12)$, 
\[\limsup_{n\to\infty}\sup_{B_{2r}}p_n|x|^{2+2\al}v_n(x)^{p_n-1}=\infty.\]Then up to a subsequence, there exists $r_n\to 0$ and $\bar z_n\in \overline{B}_{r_n}$ such that
	$$p_n|\bar z_n|^{2+2\al}v_n(\bar z_n)^{p_n-1}=\sup_{B_{r_n}}p_n|x|^{2+2\al}v_n(x)^{p_n-1}\to\iy,\quad\text{as }n\to\infty.$$
Let $\tilde z_n:=|z_n|\bar z_n$. Then
	$$\f{|\tilde z_n|}{|z_n|}=|\bar z_n|\to0  \quad\text{and}\quad  p_n|\tilde z_n|^{2+2\al}u_n(\tilde z_n)^{p_n-1}=p_n|\bar z_n|^{2+2\al}v_n(\bar z_n)^{p_n-1}\to\iy.$$
Consequently by Step 1, \[\displaystyle\limsup_{n\to\iy}\int_{B_{\delta|\tilde z_n|}(\tilde z_n)}p_n|x|^{2\al}V_n(x)u_n(x)^{p_n}\rd x\ge 4\pi e,\quad\text{for every}~\delta>0.\] Furthermore, for each fixed $\delta\in (0,1)$, we have \[B_{\delta|z_n|}(z_n)\cap B_{\delta|\tilde z_n|}(\tilde z_n)=\emptyset\quad\text{ for $n$ large},\]
so
$$\limsup_{n\to\iy}\int_{B_{\delta|z_n|}(z_n)\cup B_{\delta|\tilde z_n|}(\tilde z_n)}p_n|x|^{2\al}V_n(x)u_n(x)^{p_n}\rd x\ge 8\pi e.$$

Keep on repeating the alternatives above for the scaled new sequence (still called $v_n$) where in \eqref{al_n} we replace $z_n$ with the new sequence $\tilde z_n$, and so on. We see that, each time the scaled new sequence $v_n$ fails to verify \eqref{tem-7} for any $\e_0'\in (0,\frac12)$, we add a contribution of $4\pi e$ to the value $\int_{B_1}p_n|x|^{2\al}V_n(x)u_n(x)^{p_n}\rd x\le C$. So after finitely many steps we find a sequence $\{z_{n,1}\}$ and $\e_0'\in(0,\f{1}{2})$ such that \eqref{tem-5}-\eqref{tem-6} hold. This proves Step 2.

\vskip0.1in
{\bf Step 3.} Suppose the alternative (i) does not hold for every $\e_0\in(0,\f{1}{2})$, we prove that the alternative (ii) holds.

First, by Step 2, there are a sequence $\{z_{n,1}\}$ and $\e_0'\in(0,\f{1}{2})$ such that \eqref{tem-5}-\eqref{tem-6} hold. If there exists $\e_1'\in(0,\e_0']$ such that
	$$\sup_{\f{1}{2\e_1'}|z_{n,1}|\le |x|\le 1}p_n|x|^{2+2\al}u_n(x)^{p_n-1}\le C,$$
then by replacing $\e_0'$ with $\e_1'$, we see that the alternative (ii) holds with $l=1$. Otherwise, for any $\e\in (0,\e_0']$, 
\[\limsup_{n\to\infty}\sup_{\f{1}{2\e}|z_{n,1}|\le |x|\le 1}p_n|x|^{2+2\al}v_n(x)^{p_n-1}=\infty.\]Then up to a subsequence, there exist $\e_n\to 0$ and $y_n$ satisfying $\f{1}{2\e_n}|z_{n,1}|\le |y_n|\le 1$
such that
	$$p_n|y_n|^{2+2\al}u_n(y_n)^{p_n-1}=\sup_{\f{1}{2\e_n}|z_{n,1}|\le |x|\le 1}p_n|x|^{2+2\al}u_n(x)^{p_n-1}\to\iy,\quad\text{as}~n\to\iy.$$
This implies $\frac{|z_{n,1}|}{|y_n|}\to 0$. Furthermore, by Step 1, we see that necessarily
\be y_n\to0,\quad \limsup_{n\to\iy}\int_{B_{\delta|y_n|}(y_n)}p_n|x|^{2\al}V_n(x)u_n(x)^{p_n}\rd x\ge 4\pi e, \ee
for every $\delta>0$. To obtain the second sequence $z_{n,2}$, for $\e\in(0,\f{1}{2})$ we consider
\be\lab{tem-12}  \sup_{\lbr{ \f{1}{2\e}|z_{n,1}|\le |x|\le 2\e|y_n| }} p_n|x|^{2+2\al}u_n(x)^{p_n-1}. \ee
If there is $\e\in(0,\f{1}{2})$ such that \eqref{tem-12} is uniformly bounded for all $n$, we would simply take $z_{n,2}=y_n$, and adjust accordingly $\e_0'$ (for example, replace $\e_0'$ with $\min\{\e_0', \e\}$) in order to ensure \eqref{tem-6} with $i=1$. Otherwise, 
for any $\e\in (0,\frac12)$, 
\[\limsup_{n\to\infty}\sup_{\lbr{ \f{1}{2\e}|z_{n,1}|\le |x|\le 2\e|y_n| }} p_n|x|^{2+2\al}u_n(x)^{p_n-1}=\infty.\]
Then up to a subsequence, there are $\e_n\to 0$ and $\bar{y}_n$ satisfying
$\f{1}{2\e_n}|z_{n,1}|\le |\bar{y}_n|\le 2\e_n|y_n|$ such that
\[p_n|\bar y_n|^{2+2\al}u_n(\bar y_n)^{p_n-1}=\sup_{\f{1}{2\e_n}|z_{n,1}|\le |x|\le 2\e_n|y_n|}p_n|x|^{2+2\al}u_n(x)^{p_n-1}\to\iy.\]
Therefore, we could replace $y_n$ with this new sequence $\bar y_n$ with the properties $\f{|z_{n,1}|}{|\bar y_n|}\to0$, $\f{|\bar y_n|}{|y_n|}\to0$, 
	$$\bar y_n\to0,\quad \limsup_{n\to\iy}\int_{B_{\delta|\bar y_n|}(\bar y_n)}p_n|x|^{2\al}V_n(x)u_n(x)^{p_n}\rd x\ge 4\pi e,$$
	and consider again whether \eqref{tem-12} is uniformly bounded for some $\e\in (0,\frac12)$.
Note that, as above,  each time we admit the existence of such a new sequences, we add a contribution of $4\pi e$ to the value $\int_{B_1}p_n|x|^{2\al}V_n(x)u_n(x)^{p_n}\rd x\le C$. So by repeating the same alternatives for any such new sequence, after finitely many steps we must arrive to one for which \eqref{tem-12} is uniformly bounded on $n\in\N$ for some $\e\in(0,\f{1}{2})$. Such sequence defines $z_{n,2}$, and we can adjust $\e_0'\in(0,\f{1}{2})$ accordingly in order to guarantee \eqref{tem-6} with $i=1$
and 
\[\sup_{\f{1}{2\e_0'}|z_{n,1}|\le |x|\le 2\e_0'|z_{n,2}| } p_n|x|^{2+2\al}u_n(x)^{p_n-1}\leq C.\]

Finally we iterate the argument above by replacing $z_{n,1}$ with $z_{n,2}$. We are either able to check \eqref{tem-15}, \eqref{tem-15-0}, \eqref{tem-16} and \eqref{tem-17} for $l=2$ and so we are done, or obtain a third sequence $\{z_{n,3}\}$ for which we can verify \eqref{tem-15}, \eqref{tem-15-0} and \eqref{tem-17} for $i=1,2$. After finitely many steps we arrive to the desired conclusion, i.e. the alternative (ii) holds.
\ep

To handle the alternative (ii) in Proposition \ref{decomposition}, we need to estimate the energies in neck domains. 
Let $\e_0'\in(0,\f{1}{2})$ and $z_{n,i}$, $i=1,\cdots,l$, be given by the alternative (ii) in Proposition \ref{decomposition}. We define the subsets of $B_1$,
\be\lab{definition-PQ}\begin{aligned} 	
	Q_{n,i}&:=\lbr{x\in B_1:~\e_0'|z_{n,i}|\le|x|\le \f{1}{\e_0'}|z_{n,i}|},   \quad \text{for }1\leq i\leq l,\\
	P_{n,i}&:=\lbr{x\in B_1:~\f{1}{\e_0'}|z_{n,i-1}|\le|x|\le \e_0'|z_{n,i}|}, \quad \text{for }1\leq i\leq l+1,
\end{aligned}\ee
where we set $|z_{n,0}|=0$ and $|z_{n,l+1}|=1$. Then for any $n\ge1$, 
	$$B_1=\sbr{\bigcup_{i=1}^lQ_{n,i}}\cup\sbr{\bigcup_{j=1}^{l+1}P_{n,j}}\cup(B_1\setminus B_{\e_0'}). $$
We compute the integrals on each domain. Since \eqref{con-singular1} gives \be\label{max01}\sup_{B_1\setminus B_{\e_0'}}p_nu_n\le C,\ee we immediately obtain
\be\lab{tem-56-0}  \lim_{n\to\infty}\int_{B_1\setminus B_{\e_0'}}p_n|x|^{2\al}V_n(x)u_n(x)^{p_n-1+t}\rd x=0,\quad t=0,1,2. \ee
By \eqref{max} and \eqref{tem-15}, we have for $1\leq i\leq l$,
\[1\geq |z_{n,i}|^{\alpha_n}\geq \frac{1+o_n(1)}{u_n(z_{n,i})}\geq \frac{1}{\gamma}+o_n(1),\]
so up to a subsequence we may assume
\be\label{eq-ci} \lim_{n\to\infty} |z_{n,i}|^{-\alpha_n} =c_i \in [1, \gamma],\quad 1\leq i\leq l.\ee
Clearly
\be\label{eq-ci1} \gamma\geq c_1\geq c_2\cdots \geq c_l\geq 1. \ee

\bl\lab{energy}  Suppose the alternative (ii) in Proposition \ref{decomposition} holds, then there exists positive integers $N_i$, $i=1,\cdots,l$, such that for $t=0,1,2$,
\be\lab{energy-PQ}\begin{aligned}   
	\lim_{n\to\infty}\int_{Q_{n,i}}p_n|x|^{2\al}V_n(x)u_n(x)^{p_n-1+t}\rd x &=8\pi e^{\f{t}{2}}c_i^tN_i,\quad \forall\, 1\leq i\leq l, \\
	\lim_{n\to\infty}\int_{P_{n,i}}p_n|x|^{2\al}V_n(x)u_n(x)^{p_n-1+t}\rd x &=0,\quad \forall\, 1\leq i\leq l+1. 
\end{aligned}\ee
\el

\bp
We devide the proof into several steps.

{\bf Step 1.} we consider the integral on $P_{n,i}$ for $2\leq i\leq l+1$. Fix any $2\leq i\leq l+1$, we claim 
\be\lab{tem-50}   \sup_{P_{n,i}}p_n u_n\le C,\quad\forall~n\ge1. \ee
Assume by contradiction that \eqref{tem-50} does not hold, then up to a subsequence, we can take $y_n\in P_{n,i}$ such that
\be\lab{tem-51} 	p_nu_n(y_n)=\sup_{P_{n,i}}p_n u_n\to\iy.\ee
Let
\be  w_{n}(x):=|y_n|^{\al_n}u_n(|y_n|x),\quad\text{with}~\al_n=\f{2+2\al}{p_n-1}.\ee
Denote $D_0=\{x\in\R^2:~\f{1}{2}\le |x|\le 2\}$. It is easy to see that $-\Delta w_{n}=|x|^{2\al}V_n(|y_n|x)w_{n}^{p_n}$ in $D_0$. By \eqref{max}, \eqref{tem-15} and $\f{1}{\e_0'}|z_{n,i-1}|\le|y_n|\le \e_0'|z_{n,i}|$, we get 
	$$1\ge |y_n|^{\al_n}\ge\sbr{\f{1}{\e_0'}|z_{n,i-1}|}^{\al_n}\geq
	\sbr{\f{1}{\e_0'}}^{\al_n}\frac{1+o_n(1)}{u_n(z_{n,i-1})}
	\geq\f{1}{\ga}+o_n(1),$$
and hence
	$$\int_{D_0}p_n|x|^{2\al}V_n(|y_n|x)w_{n}(x)^{p_n}\rd x=|y_n|^{\al_n}\int_{\f{|y_n|}{2}\le |x|\le 2|y_n|}p_n|x|^{2\al}V_n(x)u_n(x)^{p_n}\rd x\le C.$$
Moreover, by \eqref{tem-16}-\eqref{tem-17}, we get 
\begin{align*}\sup_{D_0}p_n|x|^{2\alpha+2}w_{n}(x)^{p_n-1}=\sup_{\frac{|y_n|}{2}\leq |x|\leq 2|y_n|}p_n|x|^{2\alpha+2}u_{n}(x)^{p_n-1}\\
\le \sup_{\f{1}{2\e_0'}|z_{n,i-1}|\leq |x|\leq 2\e_0'|z_{n,i}|}p_n|x|^{2\alpha+2}u_{n}(x)^{p_n-1}\leq C.\end{align*}
From here, and noting that $|x|^{2\al}V_n(|y_n|x)$ is uniformly bounded  for $x\in D_0$, we can apply Theorem \ref{thm-BM} for $w_{n}$ in $D_0$, and conclude that the alternative (i) in Theorem \ref{thm-BM} holds, which implies $\sup_{|x|=1}p_nw_{n}(x)\le C$ for some $C>0$. Thus
	$$p_nu_n(y_n)=p_nw_n(\f{y_n}{|y_n|})|y_n|^{-\al_n}\le C,$$
which is  a contradiction with \eqref{tem-51}. 
So the claim \eqref{tem-50} holds. It follows that
\be\lab{tem-56}   \lim_{n\to\infty}\int_{P_{n,i}}p_n|x|^{2\al}V_n(x)u_n(x)^{p_n-1+t}\rd x=0,\quad t=0,1,2, \;\forall\, 2\leq i\leq l+1. \ee

{\bf Step 2.} we consider the integral on $P_{n,1}=\{x\in B_1:~|x|\le \e_0'|z_{n,1}|\}$. Let
\be  w_{n,1}(x):=|z_{n,1}|^{\al_n}u_n(|z_{n,1}|x),\quad\text{with}~\al_n=\f{2+2\al}{p_n-1}.\ee
Then it is easy to check that $w_{n,1}$ satisfies
\be\begin{cases}
	-\Delta w_{n,1}=|x|^{2\al}V_n(|z_{n,1}|x)w_{n,1}^{p_n},\quad \text{in}~B_{(1+\delta)\e_0'},\\
	\int_{B_{(1+\delta)\e_0'}}p_n|x|^{2\al}V_n(|z_{n,1}|x)w_{n,1}(x)^{p_n}\rd x\le C.
\end{cases}\ee
where $\delta>0$ is a small constant. 
Applying Theorem \ref{thm-BM} for $w_{n,1}$ in $B_{(1+\delta)\e_0'}$, there are two possibilities. If the alternative (i) of Theorem \ref{thm-BM} holds, then we have $\sup_{\overline B_{\e_0'}}p_nw_{n,1}\le C$, and hence 
\be\label{max02}\sup_{P_{n,1}}p_nu_n=|z_{n,1}|^{-\al_n}\sup_{\overline B_{\e_0'}}p_nw_{n,1}\le C.\ee It follows that
\be\lab{tem-57}   \lim_{n\to\infty}\int_{P_{n,1}}p_n|x|^{2\al}V_n(x)u_n(x)^{p_n-1+t}\rd x=0,\quad t=0,1,2.  \ee
If the alternative (ii) of Theorem \ref{thm-BM} holds, since \eqref{tem-16} implies
\be\lab{tem-54}  \sup_{B_{(1+\delta)\e_0'}}p_n|x|^{2+2\al} w_{n,1}(x)^{p_n-1}\leq \sup_{B_{2\e_0'|z_{n,1}|}}p_n|x|^{2+2\al} u_{n}(x)^{p_n-1}\le C, \ee
we see that $0$ is the only blow-up point of $p_nw_{n,1}$ in $B_{(1+\delta)\e_0'}$.
Moreover, by \eqref{tem-54}, one can check that the assumption \eqref{assume-1} holds for $w_{n,1}$, so we can apply Lemma \ref{tem-41} to $w_{n,1}$ to  conclude that
	$$\lim_{n\to\infty}\int_{B_{\e_0'}}p_n|x|^{2\al}V_n(|z_{n,1}|x)w_{n,1}(x)^{p_n-1+k}\rd x=8\pi(1+\al)e^{\frac{k}{2}},\quad k=0,1,2,$$
	\[\max_{B_{\e_0'}}w_{n,1}\to \sqrt{e}.\]
Backing to $u_n$ and using \eqref{eq-ci}, we get
\be\lab{tem-58} \begin{aligned}
& \lim_{n\to\infty}\int_{P_{n,1}}p_n|x|^{2\al}V_n(x)u_n(x)^{p_n-1+k}\rd x\\
=&\lim_{n\to\infty}|z_{n,1}|^{-k\alpha_n}\int_{B_{\e_0'}}p_n|x|^{2\al}V_n(|z_{n,1}|x)w_{n,1}(x)^{p_n-1+k}\rd x\\
= &8\pi (1+\alpha) c_1^k e^{\frac{k}{2}}\geq 8\pi (1+\alpha) e^{\frac{k}{2}},\quad k=0,1,2,
\end{aligned}\ee
\be\label{max03} \max_{P_{n,1}}u_n=|z_{n,1}|^{-\alpha_n}\max_{B_{\e_0'}}w_{n,1}\to c_1\sqrt{e}.\ee
We will show in Step 4 that \eqref{tem-58} can not hold, namely actually \eqref{tem-57} holds.

\vskip0.1in
{\bf Step 3.} we consider the integral on $Q_{n,i}$ for $1\leq i\leq l$. Let
\be  w_{n,i}(x):=|z_{n,i}|^{\al_n}u_n(|z_{n,i}|x),\quad\text{with}~\al_n=\f{2+2\al}{p_n-1}.\ee
By \eqref{tem-15}, it is easy to check that $w_{n,i}$ satisfies
\be\begin{cases}
	-\Delta w_{n,1}=|x|^{2\al}V_n(|z_{n,i}|x)w_{n,i}^{p_n},\quad \text{in}~D_0':=\lbr{x\in\R^2:~\e_0'\le |x|\le \f{1}{\e_0'}},\\
	\int_{D_0'}p_n|x|^{2\al}V_n(|z_{n,i}|x)w_{n,i}(x)^{p_n}\rd x\le C,\\
	p_nw_{n,i}(\f{z_{n,i}}{|z_{n,i}|})\to\iy.
\end{cases}\ee
Note from \eqref{tem-16}-\eqref{tem-17} that
\begin{align}\lab{tem-55}  &\sup_{D_0'\setminus\{2\e_0'\le |x|\le \f{1}{2\e_0'}\} }p_n|x|^{2+2\al} w_{n,i}(x)^{p_n-1}\\
=&\sup_{\{\frac{|z_{n,i}|}{2\e_0'}\leq |x|\leq \frac{|z_{n,i}|}{\e_0'}\}\cup\{\e_0'|z_{n,i}|\leq |x|\leq 2\e_0'|z_{n,i}|\}}p_n|x|^{2+2\al} u_{n}(x)^{p_n-1}\le C.\nonumber \end{align}
Therefore by applying Theorem \ref{thm-BM}, the blow-up set $\Sigma_i$ of $w_{n,i}$ in $D_0'$ has the following property:
	$$\Sigma_i\neq\emptyset\;\text{is a finite set},\quad \lim_{n\to\infty}\frac{z_{n,i}}{|z_{n,i}|}\in\Sigma_i\subset\{2\e_0'\le |x|\le \f{1}{2\e_0'}\}\Subset D_0'.$$
Hence, noting that $|x|^{2\al}V_n(|z_{n,i}|x)$ satisfies \eqref{con-regular2} for $x\in D_0'$, we are in position to apply Theorem \ref{thm-regular} to $w_{n,i}$ around each point of $\Sigma_i$ and derive
	$$\lim_{n\to\infty}\int_{D_0'}p_n|x|^{2\al}V_n(|z_{n,i}|x)w_{n,i}(x)^{p_n-1+k}\rd x=8\pi e^{\frac{k}{2}} N_i,\quad k=0,1,2,$$
	\[\max_{D_0'}w_{n,i}\to \sqrt{e},\]
where $N_i=\#\Sigma_i\geq 1$. Backing to $u_n$ and using \eqref{eq-ci}, we get
\be\label{4-56}\begin{aligned}
&\lim_{n\to\infty}\int_{Q_{n,i}}p_n|x|^{2\al}V_n(x)u_n(x)^{p_n-1+k}\rd x\\
=&\lim_{n\to\infty}|z_{n,i}|^{-k\alpha_n}\int_{D_0'}p_n|x|^{2\al}V_n(|z_{n,i}|x)w_{n,i}(x)^{p_n-1+k}\rd x\\
= &8\pi c_i^{k} e^{\frac{k}{2}} N_i\geq 8\pi e^{\frac{k}{2}} N_i,\quad k=0,1,2,
\end{aligned}\ee
\be\label{max04}\max_{Q_{n,i}}u_n=|z_{n,i}|^{-\al_n}\max_{D_0'}w_{n,i}\to c_i\sqrt{e}.\ee

{\bf Step 4.} We claim that $\gamma=c_1\sqrt{e}\geq \sqrt{e}$, and \eqref{tem-58} can not hold in Step 2, so \eqref{tem-57} holds.

Indeed, by \eqref{max01}, \eqref{tem-50}, \eqref{max02}, \eqref{max03},
\eqref{max04} and \eqref{eq-ci1}, we have
\[\gamma=\lim_{n\to\infty}\max_{B_1}u_n=\lim_{n\to\infty}\max_{P_{n,1}\cup\cup_i Q_{n,i}} u_n=\max_{1\leq i\leq l} c_i\sqrt{e}=c_1\sqrt{e}.\]
Assume by contradiction that \eqref{tem-58} hold in Step 2. Recalling $\beta_k$ defined by \eqref{tem-60-beta}, we see from \eqref{tem-56-0}, \eqref{tem-56}, and \eqref{4-56} that
\[\beta_1=8\pi (1+\alpha) c_1 e^{\frac{1}{2}}+\sum_{i=1}^l 8\pi c_i e^{\frac{1}{2}} N_i>8\pi (1+\alpha) c_1 e^{\frac{1}{2}}=8\pi (1+\alpha)\gamma,\]
a contradiction with \eqref{beta12} which says that $\beta_1\leq 8\pi(1+\alpha)\gamma$. 
Thus we finish the proof.
\ep

\bl\lab{key}
	Suppose the alternative (ii) in Proposition \ref{decomposition} holds, then there must be $l=1$, $N_1=1+\al$ and 
\be\lab{tem-202} \lim_{n\to\infty}\int_{B_1}p_n|x|^{2\al}V_n(x)u_n(x)^{p_n-1+t}\rd x =8\pi(1+\al) e^{\f{t}{2}}c_1^t,\quad t=0,1,2. \ee
In particular, $\alpha\in \mathbb{N}$ since $N_1$ is a positive integer.
\el
\bp We discuss two cases separately.

{\bf Case 1:} $l=1$. Recall the definition of $\beta_k$ in \eqref{tem-60-beta}, then \eqref{beta12} tells us that $\beta_1^2=8\pi(1+\al)\beta_2$. Since $l=1$, by \eqref{tem-56-0} and \eqref{energy-PQ}, we get
	$$\beta_k=8\pi e^{\f{k}{2}}N_1c_1^k, \quad k=0,1,2,$$
so that we obtain $c_1^2N_1^2=(1+\al)c_1^2N_1$, which gives $N_1=1+\al$ and \eqref{tem-202}.

{\bf Case 2:} $l\ge2$. We first claim that $N_1=1+\al$. We define $r_n=|z_{n,2}|$ and $v_n(x)=r_n^{\al_n}u_n(r_nx)$ with $\al_n=\f{2+2\al}{p_n-1}$. Then it is easy to see that $v_n$ satisfies
	$$-\Delta v_n=|x|^{2\al}\wt V_nv_n^{p_n},\quad\text{in}~B_{\e_0'},$$
where $\wt V_n(x)=V_n(r_nx)$. Let
\begin{align}\label{tem-60-beta2}
\wt\beta_k:=\lim_{n\to\infty} \int_{B_{\e_0'}}p_n|x|^{2\al}\wt V_nv_n(x)^{p_n-1+k}\rd x,\quad k=0,1,2. 
\end{align}
By \eqref{energy-PQ}, we get
\[\wt\beta_k =\lim_{n\to\infty}r_n^{k\al_n}\int_{P_{n,1}\cup Q_{n,1}\cup P_{n,2}}p_n|x|^{2\al}V_n(x)u_n(x)^{p_n-1+k}\rd x =8\pi e^{\frac{k}{2}} N_1c_1^{k}c_2^{-k}.\]
Since $\sup_{P_{n,2}}p_n u_n\le C$, we obtain 
	\[p_n v_n\le C, \quad\text{for any}~x\in{B_{\e_0'}\setminus B_{\f{|z_{n,1}|}{\e_0'|z_{n,2}|}}}.\]
Thanks to $\f{|z_{n,1}|}{|z_{n,2}|}\to0$, we know that $0$ is the only blow up point of $v_n$ in $B_{\e_0'}$. Now we are in the same situation of Lemma \ref{tem-203}, so that it holds $\wt\beta_1^2=8\pi(1+\al)\wt\beta_2$. It follows that $N_1=1+\al$. 

Now as in Case 1, by \eqref{tem-56-0} and \eqref{energy-PQ}, we get
	$$\beta_k=8\pi e^{\f{k}{2}}\sum_{i=1}^lN_ic_i^k,\quad k=0,1,2.$$
Then by $\beta_1^2=8\pi(1+\al)\beta_2$, we get $(\sum_{i=1}^lN_ic_i)^2=N_1\sum_{i=1}^lN_ic_i^2$. Thanks to \eqref{eq-ci1}, we have
	$$N_1\sum_{i=2}^lN_ic_i^2=(\sum_{i=1}^lN_ic_i)^2-N_1^2c_1^2>2N_1c_1\sum_{i=2}^lN_ic_i\ge 2N_1\sum_{i=2}^lN_ic_i^2, $$ 
which is a contradiction, so $l=1$ and we finish the proof. 
\ep

\bp[Proof of Theorem \ref{thm-singular}]
If $\al\not\in\N$, from Proposition \ref{decomposition} and Lemma \ref{key}, we know that \eqref{bound} holds, so that Lemma \ref{tem-41} gives the theorem. If $\al\in\N$, we don't know whether \eqref{bound} holds or not, but anyway this theorem follows from Lemma \ref{tem-41} and Lemma \ref{key}.
\ep

\section{The boundary value problem}

This section is devoted to the proof of Theorem \ref{thm-1}.
Let $u_n$ be a solution sequence of \eqref{equ-1-0}, and $W_n(x)$ satisfy \eqref{con-1}-\eqref{con-2}. We denote $\nm{u}_\iy=\nm{u}_{L^\iy(\Omega)}$ for simplicity. Let $\vp_1>0$, $\nm{\vp_1}_\iy=1$, be the eigenfunction of
the first eigenvalue of $-\Delta$ in $\Omega$ with the Dirichlet boundary condition:	$$\la_1(\Omega):=\inf_{u\in H_0^1(\Omega)}\f{\int_\Omega|\nabla u|^2}{\int_\Omega u^2}>0.$$
Then $\vp_1$ satisfies $-\Delta\vp_1=\la_1(\Omega)\vp_1$ and we have  
	$$\int_\Omega (W_nu_n^{p_n-1}-\la_1(\Omega))u_n\vp_1=\int_\Omega (-\vp_1\Delta u_n+u_n\Delta\vp_1 )=0.$$
So $(\nm{W_nu_n^{p_n-1}}_\iy-\la_1(\Omega))\int_\Omega u_n\vp_1\ge0$, which implies \[\nm{u_n}_\iy\ge \sbr{\f{\la_1(\Omega)}{\max_\Omega W_n}}^{\f{1}{p_n-1}}\geq \sbr{\f{\la_1(\Omega)}{C}}^{\f{1}{p_n-1}}.\] As a result, we obtain $\liminf_{n\to+\iy}\nm{u_n}_\iy\ge1$, which yields that $\SR\neq\emptyset$, where $\SR$ is the set of blow-up points of $p_nu_n$ defined in \eqref{S}. Moreover, the same argument as Lemma \ref{converge-regular} implies $\limsup_{n\to+\iy}\nm{u_n}_\iy\le C$.

Applying Theorem \ref{thm-BM} for $u_n$ in $\Omega$, we obtain a set 
\[\Sigma=\lbr{a_1,\cdots,a_k}\subset\Omega\] satisfying the properties in Theorem \ref{thm-BM}, where we set $\Sigma=\emptyset$ and $k=0$ if the alternative (i) holds. Obviously, $\Sigma\subset\SR$. 

Recall the zero set $\ZR\subset\Omega$ of $W_n(x)$ defined in \eqref{Z}.  We choose $r_0>0$ small such that
\be
	B_{2r_0}(a)\subset\Omega\quad\text{and}\quad B_{2r_0}(a)\cap B_{2r_0}(a')=\emptyset,\quad\text{for}~a,a'\in\Sigma\cup\ZR,~a\neq a'.
\ee
When $\Sigma\neq\emptyset$, we
define the local maximums $\ga_{n,i}$ and the local maximum points $x_{n,i}$ of $u_n$ near $a_i\in\Sigma$ by
\be\label{5--2}
	\ga_{n,i}=u_n(x_{n,i}):=\max_{B_{2r_0}(a_i)} u_n,\quad\text{for}~i=1,\cdots,k.
\ee
Then the proof of Theorem \ref{thm-BM} shows that $x_{n,i}\to a_i$ and $\ga_{n,i}\to \ga_i\geq 1$.

Now we aim to exclude the boundary concentration. Denote 
\be 	\Omega_\delta:=\lbr{x\in\Omega:~d(x,\pa\Omega)\ge\delta}. \ee
Since $\ZR\subset\Omega$ and $\Sigma\subset\Omega$, we take $\delta_0>0$ small such that $\Omega_{3\delta_0}$ is a compact subset satisfying
\be\lab{label-34} 	(\ZR\cup\Sigma)\subset\Omega_{3\delta_0}\Subset\Omega.\ee
Denote 
\[\widetilde{\Omega}:=\Omega\setminus\Omega_{2\delta_0}=\lbr{x\in\Omega:~d(x,\pa\Omega)<2\delta_0}.\] Since Theorem \ref{thm-BM} says that 
\be\label{777}\text{\it $p_nu_n$ is uniformly bounded in $L_{loc}^{\infty}(\Omega\setminus\Sigma)$},\ee we see from \eqref{label-34} that
\be\lab{label-20}  0\le p_nu_n(x)\le C,\quad\text{for}~x\in\pa\widetilde\Omega, \ee
where $C>0$ is a constant. More precisely,
\be\lab{equ-1-2} \begin{cases}
	-\Delta u_n=W_n(x)u_n^{p_n},\quad u_n>0,\quad\text{in}~\widetilde\Omega,\\
	u_n=0,\quad\text{on}~\pa\Omega\subset\pa\widetilde{\Omega},\\
	0<p_nu_n\leq C\quad\text{on}~\pa\widetilde{\Omega}\setminus\pa\Omega,\\
	\|u_n\|_{L^\infty(\Omega)}\leq C,\\
	\int_\Omega p_n W_n(x)u_n^{p_n}\rd x\le C.
\end{cases}\ee
Furthermore, it follows from \eqref{con-1}-\eqref{con-2} and \eqref{label-34} that
\be\lab{con-2-0}  	0<\f{1}{C}\le W_{n}(x)\le C<\iy,\quad |\nabla  W_{n}(x)|\le C,\quad\text{for}~x\in\widetilde\Omega.  \ee
 
\bpr\lab{boundary}
	There is no blow-up point in $\overline\Omega\setminus\Omega_{\delta_0}$ for $p_nu_n$.
\epr

To give the proof of Proposition \ref{boundary}, we assume by contradiction that there exists a blow-up point in $\overline\Omega\setminus\Omega_{\delta_0}\subset\widetilde{\Omega}$, then
\be\lab{label-21}  \limsup_{n\to\iy}\sup_{\widetilde\Omega} p_nu_n=\iy, \ee
and the same argument as \eqref{4-42} implies
\[\limsup_{n\to\iy}\sup_{\widetilde\Omega}p_nW_n(x)u_n(x)^{p_n-1}=\infty.\]
Thus up to a subsequence, there is a family of points $\{y_{n,1}\}$ such that
\be\label{con-2-22} p_nW_n(y_{n,1})u_n(y_{n,1})^{p_n-1}\to\iy,\quad\text{as}~n\to\iy.\ee

Now we suppose there exist $m\in\N\setminus\{0\}$ families of points $\{y_{n,i}\}_{n\ge1}$, $i=1,\cdots,m$, in $\widetilde\Omega$ such that
\be\lab{label-13}  p_nW_n(y_{n,i})u_n(y_{n,i})^{p_n-1}\to\iy,\quad\text{as}~n\to\iy. \ee
Define the parameters $\e_{n,i}>0$ by 
\be\lab{laben-14}  \e_{n,i}^{-2}:=p_nW_n(y_{n,i})u_n(y_{n,i})^{p_n-1},\quad\text{for}~i=1,\cdots,m, \ee
then
\be\label{5-13} \lim_{n\to\infty}\e_{n,i}=0,\quad\quad \liminf_{n\to\iy} u_n(y_{n,i})\ge 1. \ee
Define 
\be\label{5-14}  R_{n,m}(x):=\min_{i=1,\cdots,m}|x-y_{n,i}|,\quad\text{for}~x\in\widetilde\Omega. \ee
As in \cite{LE-7}, we introduce the following properties:
\begin{itemize}
\item[($\PR_1^m$)] For any $1\leq i,j\leq m$ and $i\neq j$, 	
	$$\lim_{n\to\infty}\e_{n,i}=0,\quad\quad\lim_{n\to\iy}\f{|y_{n,i}-y_{n,j}|}{\e_{n,i}}=\iy.$$
\item[($\PR_2^m$)] For each $1\leq i\leq m$, for $x\in\widetilde\Omega_{n,i}:=\f{\widetilde\Omega-y_{n,i}}{\e_{n,i}}$, 
	\be \label{5--15} w_{n,i}(x):=p_n\sbr{\f{u_n(y_{n,i}+\e_{n,i}x)}{u_n(y_{n,i})}-1}\to U_0(x)=-2\log\sbr{1+\f{1}{8}|x|^2} \ee
	in $\CR_{loc}^2(\R^2)$ as $n\to\iy$.
\item[($\PR_3^m$)] There exists $C>0$ independent of $n$ such that
	$$\sup_{x\in\widetilde{\Omega}}p_nR_{n,m}(x)^2W_n(x)u_n(x)^{p_n-1}\le C,\quad\forall\,n.$$
\end{itemize}
It is easy to see that once ($\PR_1^m$)-($\PR_3^m$) hold for $m\in\N\setminus\{0\}$ families of points $\{y_{n,i}\}_{n\ge1}$, then we can not find an $m+1$ family of points $\{y_{n,m+1}\}_{n\ge1}$ such that ($\PR_1^{m+1}$) holds.

\bl\lab{label-31}
	There exists $l\in\N\setminus\{0\}$ families of points $\{y_{n,i}\}_{n\ge1}$ in $\widetilde\Omega$, $i=1,\cdots,l$, such that, after passing to a subsequencce, $(\PR_1^l)$-$(\PR_3^l)$ hold.
\el

\bp Thanks to \eqref{equ-1-2}, \eqref{con-2-0} and \eqref{con-2-22}, the proof of Lemma \ref{label-31} is very similar to that of \cite[Proposition 2.1]{Druet} or \cite[Proposition 2.2]{LE-7}, so we omit it. 
\ep

Now we define the concentration set $\TR$ in $\overline{\widetilde{\Omega}}$ by
\be\lab{label-33}  \TR:=\lbr{a_{k+1},\cdots,a_{k+l}}=\lbr{\lim_{n\to\iy}y_{n,i},~i=1,\cdots,l}\subset\overline{\widetilde \Omega}, \ee
with $y_{n,i}\to a_{k+i}$ as $n\to\iy$, where $y_{n,i}$ is given by Lemma \ref{label-31}.

\bl\label{lem-53} We have $\TR\subset\pa\Omega$ and $p_nu_n$ is uniformly bounded in $L_{loc}^{\infty}(\overline{\widetilde{\Omega}}\setminus \TR)$.
\el

\begin{proof}
If there exists $a_{k+i}\in\TR\cap\Omega$, then by \eqref{777} and \eqref{5-13} we have $a_{k+i}\in\Sigma$. However, by the choice of $\delta_0$ in \eqref{label-34}, it holds $\TR\cap\Sigma=\emptyset$, a contradiction to $a_{k+i}\in\TR\cap\Sigma$. Thus $\TR\subset\pa\Omega$. 

Recalling \eqref{equ-1-2}, we set $p_nu_n=\phi_n+\psi_n$ with
 \[\begin{cases}
-\Delta \phi_n=p_nW_n(x)u_n^{p_n},\quad\text{in}~\widetilde\Omega,\\
\phi_n=0\quad\text{on }\;\widetilde\Omega,
\end{cases}
\quad
\begin{cases}
-\Delta \psi_n=0\quad\text{in }\; \widetilde\Omega,\\
\bar \psi=p_nu_n\in [0,C]\quad\text{on }\;\pa \widetilde\Omega.
\end{cases}
\]
Then $\|\psi_n\|_{L^\infty(\widetilde\Omega)}\leq C$. We claim that
\be\label{888}
|\nabla \phi_n(x)|\leq \frac{C}{R_{n,l}(x)},\quad \forall x\in\widetilde{\Omega}.\ee

To prove \eqref{888}, we fix any $x\in\widetilde{\Omega}$. Recalling \eqref{5-14}, we take $j_0$ such that
\[R_{n,l}(x)=\min_{i=1,\cdots,m}|x-y_{n,i}|=|x-y_{n,j_0}|.\] 
Let $\tilde{G}(x,y)$ be the Green function of $-\Delta$ in $\widetilde\Omega$ with the Dirichlet boundary condition. Then
\begin{align*}
|\nabla\phi_n(x)|&=\left|\int_{\widetilde{\Omega}}\nabla_x\widetilde{G}(x,y)p_nW_n(y)u_n(y)^{p_n}dy\right|\\
&\leq \int_{\widetilde{\Omega}\cap \{|y-x|\geq\frac{R_{n,l}(x)}{2}\}} \frac{1}{|x-y|}p_nW_n(y)u_n(y)^{p_n}dy\\
&\quad+\int_{\widetilde{\Omega}\cap \{|y-x|\leq\frac{R_{n,l}(x)}{2}\}} \frac{1}{|x-y|}p_nW_n(y)u_n(y)^{p_n}dy=I_1+I_2.
\end{align*}
By \eqref{equ-1-2} we see that $I_1\leq \frac{C}{R_{n,l}(x)}$. 

For $|y-x|\leq \frac{R_{n,l}(x)}{2}=\frac{|x-y_{n,j_0}|}{2}$, we have
\[|y-y_{n,i}|\geq |x-y_{n,i}|-|x-y|\geq \frac{R_{n,l}(x)}{2},\quad\forall i,\]
so $R_{n,l}(y)\geq \frac{R_{n,l}(x)}{2}$. Then by \eqref{equ-1-2} and ($\PR_3^l$), we have
\begin{align*}
p_nW_n(y)u_n(y)^{p_n}\leq \frac{4C}{R_{n,l}(x)^2}p_n R_{n,l}(y)^2W_n(y)u_n(y)^{p_n-1}\leq \frac{C}{R_{n,l}(x)^2},
\end{align*}
and then
\[I_2\leq \frac{C}{R_{n,l}(x)^2}\int_{|y-x|\leq\frac{R_{n,l}(x)}{2}} \frac{1}{|x-y|}dy\leq \frac{C}{R_{n,l}(x)}.\]
Therefore, \eqref{888} holds. From here and $\phi_n=0$ on $\pa\widetilde\Omega$, we see that $\phi_n$ is uniformly bounded in $L_{loc}^{\infty}(\overline{\widetilde{\Omega}}\setminus \TR)$ and so does $p_nu_n=\phi_n+\psi_n$.
\end{proof}


Recalling the set $\SR$ of blow-up points defined by \eqref{S}, clearly we have
\be\lab{label-37} \SR=\Sigma\cup\TR=\lbr{a_1,\cdots,a_{k+l}}. \ee
Recalling \eqref{5--2}, \eqref{laben-14}, \eqref{5--15} and Lemma \ref{label-31}, we unify the notations by setting
	$$x_{n,k+i}:=y_{n,i},\quad \mu_{n,k+i}:=\e_{n,i},\quad v_{n,k+i}:=w_{n,i},\quad \ga_{n,k+i}:=u_n(y_{n,i})$$
for $1\leq i\leq l$, and after passing to a subsequence, we assume
	$$\ga_{n,i}\to \ga_i\ge1,\quad\text{as}~n\to\iy, \quad\text{for }1\leq i\leq k+l.$$
By Theorem \ref{thm-BM} and Lemma \ref{lem-53}, we conclude that
\be\label{5-19}\nm{p_nu_n}_{L^\iy(K)}\le C_K\quad\text{\it for any compact subset $K\in\overline\Omega\setminus\SR$}.\ee
More precisely, we have the following convergence result.
\bl\lab{converge-2}
	There exists $\sigma_{a_i}\geq 8\pi$, $i=1,\cdots,k+l$, such that
	$$p_nu_n(x)\to \sum_{i=1}^{k+l}\sigma_{a_i} G(x,a_i),\quad\text{in}~\CR_{loc}^2(\overline\Omega\setminus\SR)~\text{as}~n\to\iy. $$
\el
\bp
Exactly as in Lemma \ref{converge-regular2}, we get
	$$p_nu_n(x)\to \sum_{i=1}^{k+l} \sigma_{a_i} G(x,a_i),\quad\text{in}~\CR_{loc}^2(\overline{\Omega}\setminus \SR),\quad\text{as}~n\to\iy,$$
where 
	$$\sigma_{a_i}:=\lim_{d\to0}\lim_{n\to\iy}p_n\int_{B_d(a_i)\cap\Omega}W_n(x)u_n(x)^{p_n}\rd x.$$
We show that $\sigma_{a_i}\geq 8\pi$. For $1\leq i\leq k$, by Theorems \ref{thm-regular} and \ref{thm-singular}, we see that $\sigma_{a_i}\geq 8\pi \sqrt{e}$. For $k+1\leq i\leq k+l$, since $x_{n,i}\to a_i$ as $n\to\iy$, then $B_{d/2}(x_{n,i})\subset B_d(a_i)$ for $n$ large, and hence
$$\begin{aligned}
	p_n\int_{B_d(a_i)\cap\Omega}W_n(x)u_n(x)^{p_n}\rd x
	&\ge p_n\int_{B_{d/2}(x_{n,i})\cap\Omega}W_n(x)u_n(x)^{p_n}\rd x\\
	&=\ga_{n,i}\int_{B_{\f{d}{2\mu_{n,i}}}\cap\Omega_{n,i}}\f{W_n(x_{n,i}+\mu_{n,i}x)}{W_n(x_{n,i})}\sbr{1+\f{v_{n,i}}{p_n}}^{p_n}\rd x.
\end{aligned}$$
Passing to the limit as $n\to\iy$, since $B_{\f{d}{2\mu_{n,i}}}\cap\Omega_{n,i}\to\R^2$, thanks to $(\PR_2^l)$ we get
\be\lab{tem-100}  \lim_{n\to\iy}p_n\int_{B_d(a_i)\cap\Omega}W_n(x)u_n(x)^{p_n}\rd x\ge \ga_i\int_{\R^2}e^{U_0}\rd x\ge 8\pi\ga_i,\ee
which gives $\sigma_{a_i}\ge 8\pi$.
\ep

Now we are ready to show that there is no boundary blow up.
\bp[Proof of Proposition \ref{boundary}]
Assume by contradiction that there exists a concentration point in $\overline\Omega\setminus\Omega_{\delta_0}$. Then the above argument shows that the blow-up point $\SR$ is given in \eqref{label-37}. Since $\TR\neq\emptyset$ now, we take $a_i\in\TR\subset\pa\Omega$ for some $k+1\leq i\leq k+l$. Choose $r>0$ such that $\SR\cap B_r(a_i)=\{a_i\}$. 

Let $y_n=a_i+\rho_{n,d}\nu(a_i)$, where $\nu(x)$ is the outer normal vector of $\pa\Omega$ at the point $x\in\pa\Omega$, and
	$$\rho_{n,d}:=\f{ \int_{\pa\Omega\cap B_d(a_i)}\sbr{\f{\pa u_n(x)}{\pa\nu}}^2\abr{x-a_i,\nu(x)}\rd s_x } { \int_{\pa\Omega\cap B_d(a_i)}\sbr{\f{\pa u_n(x)}{\pa\nu}}^2\abr{\nu(a_i),\nu(x)}\rd s_x  }.$$
Recall the zero set $\ZR$ of $W_n(x)$ defined in \eqref{Z}.
Choose $0<d<\min\{r,\f{d(a_i,\ZR)}{2}\}$ small enough such that $\f{1}{2}\le \abr{\nu(a_i),\nu(x)}\le 1$ for $x\in \pa\Omega\cap B_d(a_i)$. With this choice of $d$, we have 
\be\lab{label-40}   |\rho_{n,d}|\le 2d. \ee
Moreover, it is easy to see that the choice of $y_n$ implies
\be\lab{label-39} \int_{\pa\Omega\cap B_d(a_i)}\sbr{\f{\pa u_n(x)}{\pa\nu}}^2\abr{x-y_n,\nu(x)}\rd s_x=0.\ee
Applying the local Pohozaev identity \eqref{pho-1} in the set $\Omega\cap B_d(a_i)$ with $u=u_n$, $V(x)=W_n(x)$ and $y=y_n$, using \eqref{label-39}, the boundary condition $u_n=0$ on $\pa\Omega$ (so that $\nabla u_n=-|\nabla u_n|\nu$ on $\pa\Omega$) we obtain
{\allowdisplaybreaks
\begin{align}\lab{label-41} 
	&\quad  \f{2p_n^2}{p_n+1}\int_{\Omega\cap B_d(a_i)}W_n(x)u_n(x)^{p_n+1}\rd x\nonumber\\
	&\quad+\f{p_n^2}{p_n+1}\int_{\Omega\cap B_d(a_i)}\abr{\nabla W_n(x),x-y_n} u_n(x)^{p_n+1}\rd x\nonumber\\
	&=\int_{\Omega\cap\pa B_d(a_i)} \abr{p_n\nabla u_n,\nu}\abr{p_n\nabla u_n,x-y_n} \rd s_x\\
	&\quad-\f{1}{2}\int_{\Omega\cap\pa B_d(a_i)} \abs{p_n\nabla u_n}^2\abr{x-y_n,\nu} \rd s_x \nonumber\\
	&\quad +\f{p_n^2}{p_n+1}\int_{\Omega\cap\pa B_d(a_i)} W_n(x)u_n(x)^{p_n+1}\abr{x-y_n,\nu} \rd s_x. \nonumber
\end{align} 
}%
Next we estimate the second term in the left-hand side  and all the three terms in the right-hand side.

By the choice of $d$ and \eqref{label-40}, we have $W_n(x)\ge C>0$, $|\nabla W_n(x)|\le C$ and $|x-y_n|\le 3d$ for $x\in\Omega\cap B_d(a_i)$, so
	$$\begin{aligned}
	&\abs{\f{p_n^2}{p_n+1}\int_{\Omega\cap B_d(a_i)}\abr{\nabla W_n(x),x-y_n} u_n(x)^{p_n+1}\rd x}\\
	\le & Cd\int_{\Omega\cap B_d(a_i)}p_nW_n(x)u_n(x)^{p_n}\rd x=O(d).
	\end{aligned}$$
By \eqref{5-19} we have	
	$$\f{p_n^2}{p_n+1}\abs{\int_{\Omega\cap\pa B_d(a_i)} W_n(x)u_n(x)^{p_n+1}\abr{x-y_n,\nu} \rd s_x}\le  \f{C_d^{p_n+1}}{p_n^{p_n}}\to 0\;\text{as }n\to\infty.$$
	By \eqref{label-40} we may assume $y_n\to y_d$ with $|y_d-a_i|\leq 2d$. Recall Lemma \ref{converge-2} that \[p_nu_n(x)\to F(x):=\sum_{j=1}^{k+l}\sigma_{a_j}G(x,a_j),\quad\text{in }\CR_{loc}^2(\overline\Omega\cap B_r(a_i)\setminus\{a_i\}).\] 
Since $a_i\in\partial\Omega$ implies $G(x,a_i)\equiv 0$ for $x\neq a_i$, we have (see e.g. \cite[(3.7)]{LE-1}) 
\[F(x)=O(1),\quad \nabla F(x)=O(1),\quad \text{\it for }x\in \overline{\Omega}\cap B_{r}(a_i)\setminus\{a_i\}.\]
From here and $0<d<r$, we obtain
	\begin{align*}&\lim_{n\to\infty}\int_{\Omega\cap\pa B_d(a_i)} \abr{p_n\nabla u_n,\nu}\abr{p_n\nabla u_n,x-y_n} \rd s_x\\
	=&\int_{\Omega\cap\pa B_d(a_i)} \abr{\nabla F,\nu}\abr{\nabla F,x-y_d} \rd s_x= O(1) \int_{\Omega\cap\pa B_d(a_i)}\abs{x-y_d} \rd s_x=O(d^2),\end{align*}
and similarly	$$\lim_{n\to\infty}\int_{\Omega\cap\pa B_d(a_i)} \abs{p_n\nabla u_n}^2\abr{x-y_n,\nu} \rd s_x=O(d^2).$$
	
Inserting these estimates into \eqref{label-41}, we finally obtain
	$$\lim_{d\to0}\lim_{n\to\iy}p_n\int_{\Omega\cap B_d(a_i)}W_n(x)u_n(x)^{p_n+1}\rd x=0.$$
However, a similar argument as \eqref{tem-100} leads to
	$$\lim_{d\to0}\lim_{n\to\iy}p_n\int_{\Omega\cap B_d(a_i)}W_n(x)u_n(x)^{p_n+1}\rd x\ge 8\pi\ga_i^2>0,$$
which is a contradiction. This completes the proof. 
\ep

\bp[Proof of Theorem \ref{thm-1}]
Since Proposition \ref{boundary} tells us that $\SR=\Sigma=\{a_1,\cdots,a_k\}\subset \Omega$, by using Theorem \ref{thm-regular} for those $a_i\in \SR\setminus\ZR$ and Theorem \ref{thm-singular} for those $a_i\in \SR\cap\ZR$, one can easily prove Theorem \ref{thm-1}.
\ep

\section{The ground state of the H\'enon equation}

This section is devoted to the proof of Theorem \ref{thm-0}.
Let $u_n$ be a ground state of the H\'enon equation \eqref{equ-henon}. Set 
	$$u_n(x_n)=\max_{\Omega}u_n,$$
then \eqref{tem-4-00} implies $u_n(x_n)\to\ga\in [1,\sqrt{e}]$, i.e. $p_nu_n(x_n)\to\infty$. Applying Theorem \ref{thm-1}, we see the existence of $k\in\mathbb{N}\setminus\{0\}$ and a set $\SR=\lbr{a_1,\cdots,a_k}\subset\Omega$ consisting of blow-up points of $p_nu_n$ in $\overline{\Omega}$ such that $\max_{B_r(a_i)}u_n\to\gamma_i\geq \sqrt{e}$ for any small $r$ and $\nm{p_nu_n}_{L^\iy(K)}\le C_K$ for any compact subsets $K\subset\overline\Omega\setminus\SR$. In particular,  $u_n(x_n)\to \sqrt{e}$.
\bl\label{lem-6-1}
	It holds $\SR=\{a\}$ with $a=a_1\neq 0$. Consequently, $x_n\to a$ and
	\begin{equation}\label{6--10} p_n|x|^{2\al}u_n^{p_n-1+k}\to8\pi e^{\frac{k}{2}}\delta_{a},\quad k=0,1,2 \end{equation}
	weakly in the sense of measures.
\el
\bp
Assume by contradiction that $0\in\SR$. By choosing $r>0$ small, we know that $0$ is the only blow-up point of $p_nu_n$ in $B_{r}$, i.e.,
	$$\max_{B_r}p_nu_n\to\iy  \quad \text{and}\quad \max_{\overline B_r\setminus B_\delta} p_nu_n\le C_\delta,\quad\text{for any}~0<\delta<r.$$
Applying Theorem \ref{thm-singular}, up to a subsequence we obtain 
$$\int_{B_r} p_n|x|^{2\al}u_n^{p_n+1}\rd x\to 8\pi(1+\al)ec^2\geq 8\pi(1+\al)e,$$
a contradiction with \eqref{tem-4}.

This proves $0\not\in\SR$.
Then we can apply Theorem \ref{thm-regular} around each point of $\SR$ and obtain (note \eqref{tem-4})
\[8\pi e=\lim_{n\to\infty}\int_{\Omega}p_n|x|^{2\al}u_n^{p_n+1}\rd x=8\pi e\times \#\SR.\]
This implies $\#\SR=1$, i.e. $\SR=\{a\}$ with $a=a_1$ and so $x_n\to a$. Consequently, \eqref{6--10} follow from Theorem \ref{thm-regular}.
\ep

We need the local Pohozaev identity.
\bl
	Suppose $u$ satisfies
	$$\begin{cases} -\Delta u=V(x)u^p,\quad\text{in}~\Omega,\\ u>0,\quad\text{in}~\Omega,\end{cases}$$
	then for any $y\in\R^2$ and any subset $\Omega'\subset\Omega$, it holds
	\be\lab{pho-2} \begin{aligned} 
		&\quad \f{1}{p+1}\int_{\Omega'}\pa_iV(x)u(x)^{p+1}\rd x-\f{1}{p+1}\int_{\pa\Omega'}V(x)u(x)^{p+1}\nu_i(x)\rd s_x \\
		&=\int_{\pa\Omega'}\abr{\nabla u(x),\nu(x)}\pa_iu(x)-\f{1}{2}|\nabla u(x)|^2\nu_i(x)\rd s_x,\quad i=1,2,
	\end{aligned} \ee
	where $\partial_i=\frac{\partial}{\partial x_i}$ and $\nu(x)=(\nu_1(x),\nu_2(x))$ is the outer normal vector of $\partial \Omega'$ at $x$.
\el
\bp
By direct computations, for $i=1,2$, we have
	$$-\Delta u(x)\cdot\pa_iu(x)=-\operatorname{div}(\pa_iu(x)\nabla u(x) )+\f{\pa_i|\nabla u(x)|^2}{2},$$
and
	$$V(x)u(x)^p\cdot\pa_iu(x)=\pa_i\sbr{\f{V(x)u(x)^{p+1}}{p+1}}-\f{\pa_iV(x)u(x)^{p+1}}{p+1}.$$
Then by multiplying $-\Delta u=V(x)u^p$ with $\pa_iu(x)$, integrating on $\Omega'$ and using the divergence theorem, we obtain \eqref{pho-2}.
\ep

Now we can finish the proof of Theorem \ref{thm-0}.
\bp[Proof of Theorem \ref{thm-0}]
Thanks to Lemma \ref{lem-6-1}, we can apply Lemma \ref{converge-2} and Lemma \ref{lem-6-1} to obtain
	\be\label{6-66}p_nu_n\to 8\pi\sqrt e G(x,a),\quad\text{in}~\CR_{loc}^2(\overline{\Omega}\setminus\{a\})~\text{as}~n\to\iy. \ee
It remains to compute the location of the blow-up point $a$. Applying the Pohozaev identity \eqref{pho-2} with $y=0$, $\Omega'=B_d(a)$, $V=|x|^{2\al}$ and $u=u_n$, and by using $\max_{\pa B_d(a)}p_nu_n\le C_d$, we obtain
\be\begin{aligned}\lab{tem-09}
	&\int_{\pa B_d(a)}\abr{p_n\nabla u_n,\nu}p_n\pa_i u_n-\f{1}{2}\abs{p_n\nabla u_n}^2\nu_i\rd \sigma_x\\
	=&\f{2\al p_n^2}{p_n+1}\int_{B_d(a)}|x|^{2\al-2}x_iu_n^{p_n+1}\rd x+o_n(1),\quad i=1,2.
\end{aligned}\ee
Note from \eqref{6-66} that on $\pa B_d(a)$, 
	$$\abr{p_n\nabla u_n,\nu}p_n\pa_i u_n-\f{1}{2}\abs{p_n\nabla u_n}^2\nu_i\to-64\pi^2e\sbr{\f{(x-a)_i}{8\pi^2d^3}+\f{1}{2\pi r}\pa_iH(x,a)+O(1)},$$
as $n\to\iy$. This means that
 	$$\text{LHS of \eqref{tem-09}}=-\f{64\pi^2e}{2\pi d}\int_{\pa B_d(a)}\pa_iH(x,a)\rd \sigma_x+O(r)+o_n(1).$$
On the other hand, recalling Remark \ref{decay-regular2}, we use the Domainted Convergence Theorem to get (write $a=(a_1,a_2)$)
	$$\begin{aligned}
	&\quad\text{RHS of \eqref{tem-09}}\\
	&=\f{2\al p_n}{p_n+1}\f{u_n(x_n)^2}{|x_n|^{2\al}}\int_{B_{\f{d}{2\mu_n}}}\sbr{|x_n+\mu_ny|^{2\al-2}}(x_n+\mu_ny)_i\sbr{1+\f{v_n}{p_n}}^{p_n+1}\rd y+o_n(1)\\
	&=2\al\f{e}{|a|^{2\al}}\int_{\R^2}|a|^{2\al-2}a_ie^U\rd y+o_n(1)=16\pi e\al\f{a_i}{|a|^2}+o_n(1).
	\end{aligned}$$
Thus by letting $n\to\iy$ first and then $d\to0$ in \eqref{tem-09}, we obtain $\pa_iH(a,a)=\f{\al a_i}{4\pi|a|^2}$ for $i=1,2$, which implies
$$\nabla \sbr{R(\cdot)-\f{1}{4\pi}\log|\cdot|^{2\al}}(a)=0.$$
This completes the proof.\ep


\subsection*{Acknowledgements}  Z.Chen is supported by NSFC (No. 12222109, 12071240), and H.Li is supposed by the postdoctoral fundation of BIMSA.

\end{document}